\documentclass[12pt,reqno]{amsart}
\usepackage{amsmath,amssymb}
\usepackage{amssymb}
\usepackage{amsthm}
\usepackage{color}
\usepackage[numbers,sort&compress]{natbib}
\usepackage[colorlinks,citecolor=red,linkcolor=blue]{hyperref}
\usepackage{lineno}
\usepackage{times}
\usepackage{datetime}
\usepackage{scrtime}
\usepackage{pgfplots}
\usepackage{tikz}
\usetikzlibrary{patterns,arrows,positioning,calc,fadings,shapes,decorations.markings}
\usepackage[left=2.5cm, right=2.5cm, top=2cm, bottom=2cm]{geometry}
\begin{document}
\voffset=-1 true cm \setcounter{page} {1}
\numberwithin{equation}{section}
\newtheorem{Theorem}{\bf Theorem}[section]
\newtheorem{Lem}{\bf Lemma}[section]
\newtheorem{Prop}[Lem]{\bf Proposition}
\newtheorem{Def}[Lem]{\bf Definition}
\newtheorem{Coro}[Lem]{\bf Corollary}
\newtheorem{Rem}{\bf Remark}[section]
\newcommand{\RN}{\mathbb{R}^N}
\newcommand{\RR}{\mathbb{R}^N}
\newcommand{\Proof}{\noindent\emph{Proof. }}
\newcommand{\2}{\frac{N+\alpha}{N}}
\newcommand{\II}{\int_{\RN}(I_{\alpha}*|u|^{\2})|u|^{\2}dx}
\newcommand{\IIN}{\int_{\RN}(I_{\alpha}*|u_n|^{\2})|u_n|^{\2}dx}
\newcommand{\E}{\mathcal{E}}
\title[lower critical Choquard equation]{Normalized solutions to lower critical Choquard equation in mass-supercritical setting}
\date{}
\author[Shuai Mo]{Shuai Mo}
\address[Shuai Mo]
	{\newline\indent
			School of Mathematical Sciences and LPMC,
			\newline\indent
			Nankai University
			\newline\indent
			Tianjin 300071, China }
	\email{moshuai2021@126.com}
\author[Shiwang Ma]{Shiwang Ma$^*$}
\address[Shiwang Ma]
	{\newline\indent
			School of Mathematical Sciences and LPMC,
			\newline\indent
			Nankai University
			\newline\indent
			Tianjin 300071, China }
	\email{shiwangm@nankai.edu.cn}
\maketitle
\let\thefootnote\relax
\footnotetext{The research is supported by the National Natural Science Foundation of China (Grant Nos.11571187, 11771182).

* The Corresponding author is Shiwang Ma. }
\begin{abstract}
We study the normalized solutions $(u,\lambda) \in H^1(\RN) \times \mathbb{R}^+$ to the following Choquard equation
  \begin{equation*}
     \aligned  &-\Delta u + \lambda u =\mu g(u) + \gamma (I_\alpha * |u|^{\frac{N+\alpha}{N}})|u|^{\frac{N+\alpha}{N}-2}u   & \text{in\ \ }  \mathbb{R}^N \endaligned
    \end{equation*}
under the $L^2$-norm constraint $\int_{\mathbb{R}^N}u^2 dx =c^2$.
Here $\gamma>0$, $ N\geq 1$, $\alpha\in(0,N)$, $I_{\alpha}$ is the Riesz potential,
and the unknown parameter $\lambda$ appears as a Lagrange multiplier. In a mass supercritical setting on $g$, we find regions in the $(c,\mu)$--parameter space such that the corresponding equation admits a positive radial ground state solution.
To overcome the lack of compactness resulting from the nonlocal term, we present a novel compactness lemma and some prior energy estimate.
These results
are even new for the power type nonlinearity $g(u)= |u|^{q-2}u$ with $2+\frac{4}{N}<q<2^*$ ($2^*:=\frac{2N}{N-2}$, if $N\geq 3$ and $2^* = \infty$, if $N=1, 2$).
We also show that as $\mu$ or $c$ tends to $0$ (resp. $\mu$ or $c$ tends to $+\infty$), after a suitable rescaling the ground state solutions converge in $H^1(\RN)$
to a particular solution of the limit equations.

Further, we study the non-existence and multiplicity of positive radial solutions to
\begin{equation*}
  -\Delta u + u = \eta |u|^{q-2}u + (I_\alpha * |u|^{\frac{N+\alpha}{N}})|u|^{\frac{N+\alpha}{N}-2}u, \quad \text{in}\ \  \RN
\end{equation*}
where $N \geq 1$,  $ 2< q<2^*$ and $\eta>0$.
Based on some analytical ideas the limit behaviors of the normalized solutions, we verify some threshold regions of $\eta$ such that the
corresponding equation has no positive least action solution or admits multiple positive solutions.
To the best of our knowledge, this seems to be the first result concerning the non-existence and multiplicity of positive
solutions to Choquard type equations involving the lower critical exponent.
\end{abstract}



\section{Introduction}
We study the following mass-constraint Choquard equation
\begin{equation}\label{ThePro}
    \left\{ \aligned  &-\Delta u + \lambda u =\mu g(u) +\gamma (I_\alpha * |u|^{\frac{N+\alpha}{N}})|u|^{\frac{N+\alpha}{N}-2}u   & \text{in\ \ }  \mathbb{R}^N,\\
    & \int_{\mathbb{R}^N}u^2 dx =c^2, \quad u \in H^1(\mathbb{R}^N). \endaligned \right.
    \end{equation}
Here $N \geq 1$, $c>0$ is assigned, $g$ is a continuous function, and $\lambda$ will appear as a Lagrange multiplier. The Riesz potential $I_\alpha$ with order $\alpha \in (0,N)$ is defined for $x \in \mathbb{R}\setminus \{ 0 \}$ by
$$  I_{\alpha}(x):=\frac{A_\alpha(N)}{|x|^{N-\alpha}},\quad A_{\alpha}(N):=\frac{\Gamma(\frac{N-\alpha}{2})}{\Gamma(\frac{\alpha}{2})\pi^{N/2}2^{\alpha}}  $$
with $\Gamma$ denoting the Gamma function and $*$ denoting the convolution on $\RN$.

An important motivation driving the study of (\ref{ThePro}) is searching for standing waves
 \begin{equation}\label{formStandingWave}
   \Psi(t,x)=e^{i\lambda t}u(x)
 \end{equation}
of the nonlinear Choquard equation
\begin{equation}\label{Motivationequ}
  \left\{ \aligned &i \Psi_t + \Delta \Psi + \gamma (I_{\alpha}*|\Psi|^{p})|\Psi|^{p-2}\Psi + H (|\Psi|^2)\Psi = 0,  & (x,t)\in \RN \times \mathbb{R},\\
                   &\Psi =\Psi(x,t)\in \mathbb{C},\ N\geq 1,  \endaligned \right.
 \end{equation}
since the ansatz (\ref{formStandingWave}) leads (\ref{Motivationequ})  to
$$ -\Delta u + \lambda u =\mu g(u) +\gamma (I_\alpha * |u|^{p})|u|^{p-2}u   \quad \text{in\ \ }  \mathbb{R}^N. $$
In the case $\mu=0$ and $\lambda>0$, this equation reduces to the well-known Choquard equation
\begin{equation}\label{TheChoquard}
 \aligned & -\Delta u + \lambda u = (I_\alpha * |u|^{p})|u|^{p-2}u & \text{in\ \ }  \mathbb{R}^N.\endaligned
\end{equation}
When $N=3$, $\alpha=2$, $p=2$, Equ. (\ref{TheChoquard}) appeared in the work by S.I.Pekar\cite{Pekar} describing the model of a polaron at rest in the quantum physics.
In 1976 P. Choquard used it to give rise to the theory of one-component plasma in describing the electron trapped in its own hole \cite{Lieb1976}.
It is also known as the Schr\"{o}dinger-Newton equation or the Hartree equation in the description of self-gravitating matter
in R.Penrose \cite{Moroz1998,Penrose1,Penrose2} and in M.Lewin\cite{Lewin}.
The study of (\ref{TheChoquard}) using variational methods can be traced back to the works of E.H. Lieb \cite{Lieb1976}, P.-L. Lions \cite{Lions1980,Lions1984}, and V. Moroz, J. Van Schaftingen \cite{MorozJFA}, and now the related literature is huge, for example \cite{zhangjianjun,Tangxianhua2,Xiajiankang,Liaofangfang,Ma,MaMoroz,MorozGuide,MorozCCM,LiMa} and the references therein.
We also refer to \cite{Giulini,MorozGuide} for more physical backgrounds.

Furthermore, since $\|\Psi(t,x) \|_{L^2(\RR)}$ is preserved along the time evolution:
$$  \|\Psi(\cdot,t) \|_2^2 = \| \Psi(\cdot,0)\|_2^2 \quad \text{\ for\ any\ }t \in \mathbb{R},$$
it is particularly relevant finding solutions with a prescribed $L^2$-norm, which are often referred to as \emph{normalized solutions} or \emph{$L^2$-mass constraint problem} in the literature. Indeed, the $L^2$-norm usually carries physical meaning in various nonlinear fields such as Bose-Einstein condensate \cite{Baoweizhu} or ergodic Mean Field Games systems \cite{Postoia}.
In this direction, let $G(t)=\int_{0}^{t}g(s)ds$, and define on $H^1(\RN)$ the energy functional
\begin{equation}\label{energyfunctional}
 E_{\mu,c}(u)=\frac{1}{2}\| \nabla u \|_2^2 - \mu \int_{\RN}G(u) dx - \frac{ \gamma  }{2p} \int_{\RN}(I_{\alpha} * |u|^{p})|u|^{p}dx.
\end{equation}
It is standard to check that $E_{\mu,c} \in C^1$ under some mild assumptions on $g$, and the critical points of $E_{\mu,c}$
constrainted on the sphere
$$ S(c):=\left\{ u\in H^1(\RN): \int_{\RN}|u|^2dx =c^2 \right\}$$
give rise to a solution with a prescribed $L^2$-norm. Consequently, the parameter $\lambda$ appears as a Lagrange multiplier, and we say a couple $(u,\lambda) \in H^1(\RN) \times \mathbb{R}$ is a normalized solution.

The normalized solutions to the Choquard equation (\ref{TheChoquard}) have been studied by a plenty of literatures during last decades, we refer the readers to see \cite{TanakaCV,TanakaIndAM,YeTopo,DengJMP,GuoJMP} for the case $\frac{N+\alpha}{N} < p \leq \frac{N+\alpha+2}{N}$, and see \cite{YuanERA,XvMa,XiaZhang,SPDEA,LiYeJMP} for the case $\frac{N+\alpha+2}{N} < p <\frac{N+\alpha}{N-2}$. Besides, V. Moroz, J. Van Schaftingen \cite{MorozJFA} obtained the nonexistence result of (\ref{TheChoquard}) under the range
$$p \leq \frac{N+\alpha}{N} \quad \text{or}\quad p \geq \frac{N+\alpha}{N-2}. $$
We note that the \emph{upper critical exponent} $\frac{N+\alpha}{N-2}$ comes from the Hardy-Littlewood-Sobolev inequality, and plays a similar role as the Sobolev critical exponent in the local semilinear equations, see \cite{BrezisNirenberg}.
However, the \emph{lower critical exponent} $\frac{N+\alpha}{N}$ seems to be a new feature for Choquard’s equation \cite{MorozCCM}, which is related to a new phenomenon of ``bubbling at infinity".
So it is interesting and challenging to study (\ref{ThePro}), that is (\ref{TheChoquard}) with $p=\frac{N+\alpha}{N}$ or $p=\frac{N+\alpha}{N-2}$ perturbed by a local term $g(u)$.
In the recent papers \cite{Li1,Li2}, X. Li
considered the Choquard equation with upper critical exponent and homogeneous local perturbation
\begin{equation}\label{LiPro}
    \left\{ \aligned  &-\Delta u + \lambda u =\mu |u|^{q-2}u + (I_\alpha * |u|^{\frac{N+\alpha}{N-2}})|u|^{\frac{N+\alpha}{N-2}-2}u   & \text{in\ \ }  \mathbb{R}^N,\\
    & \int_{\mathbb{R}^N}u^2 dx =c^2, \quad u \in H^1(\mathbb{R}^N), \endaligned \right.
\end{equation}
where $N \geq 3$, $2<q<2^*$, $\mu >0$. By employing strategy utilized in the study of the Sobolev critical Schr\"{o}dinger equation \cite{JeanjeanMA,SoaveJFA,Wuyuanze}, she obtained (\ref{LiPro}) admits one radial solution when $2+\frac{4}{N} \leq  q<2^*$ and two radial solutions when $2<q<2+\frac{4}{N}$. Moreover, qualitative properties and stability of solutions are also described.

It is worth to note that the aforementioned literatures encounter a common and substantial difficulty arising from the non-compact embedding
$$ H^1(\RN) \hookrightarrow L^2(\RN). $$
Specifically, the weak limit of a Palais-Smale sequence (or minimizing sequence) does not necessarily lie on the prescribed $L^2$-sphere.
Meanwhile, it also prevents the Hardy-Littlewood-Sobolev inequality from guaranteeing
the convergence
$$ \IIN \rightarrow \II $$
if $u_n \rightharpoonup u $ in $H^1(\RN)$ even in the radially symmetry function space $H^1_{rad}(\RN)$.
This represents a new difficulty caused by the intricate interplay among the $L^2$-constraint and the lower critical term.
In light of this issue, investigating (\ref{ThePro}) will be highly challenging, but it will also yield novel and deeper insights into the $L^2$-norm constraint problem.
However, to the best of our knowledge, there are only few papers concerning the lower critical Choquard equations.

Recently, S. Yao et al. \cite{Yao}  considered the model
\begin{equation*}
    \left\{ \aligned  &-\Delta u + \lambda u =\mu |u|^{q-2}u +\gamma (I_\alpha * |u|^{\frac{N+\alpha}{N}})|u|^{\frac{N+\alpha}{N}-2}u   & \text{in\ \ }  \mathbb{R}^N,\\
    & \int_{\mathbb{R}^N}u^2 dx =c^2, \quad u \in H^1(\mathbb{R}^N). \endaligned \right. \eqno{\left(\mathcal{A}_{\mu,c}\right)}
\end{equation*}
Let $N \geq 2$, $2<q<2+\frac{4}{N}$, $\gamma>0$ and $c>0$, by solving the minimization problem
\begin{equation}\label{minimizingproblrm}
   \sigma(c)=\inf_{S(c)}E_{\mu,c}(u),
\end{equation}
they concluded there exists $\mu_0>0$ such that $(\mathcal{A}_{\mu,c})$ admits a normalized ground state solution $u_0$ if $\mu > \mu_0$ with some $\lambda_0>0$. For $q=2+\frac{4}{N}$ and $\mu$ small, they also deduced the nonexistence of normalized  ground state solution. Here the ``normalized  ground state solution'' is defined as follows:
\begin{Def}
We say that $(u,\lambda)$ is a normalized  ground state solution to (\ref{ThePro}) if
$$(E_{\mu,c}|_{S(c)})'(u)=0,\ \  E_{\mu,c}(u):=\inf\{  E_{\mu,c}(v): v\in S(c), (E_{\mu,c}|_{S(c)})'(v)=0    \}.    $$
\end{Def}
More precisely, they proved the minimizing sequence of $\sigma(c)$ converges strongly in $H^1(\RN)$ if
\begin{equation}\label{Theminiestimate}
   \sigma(c) < -\frac{\gamma N}{2(N+\alpha)} S_{\alpha}^{-\2}c^{\frac{2(N+\alpha)}{N}}
\end{equation}
where \begin{equation}\label{TheLowerConstant}
        S_{\alpha} = \inf_{u \in H^1(\RN)} \frac{\| u\|_2^2}{\left( \II  \right)^{\frac{N}{N+\alpha}}}.
      \end{equation}
It is well-known \cite[Theorem 4.3]{Lieb2001}  that  the family of extremal functions of $S_{\alpha}$ is given by
\begin{equation}\label{extrmalLower}
  {V}_{\delta}(x)=a \left( \frac{\delta}{\delta^2 + |x-y|^{2}}   \right)^{\frac{N}{2}}, \quad \delta>0,\quad y \in\mathbb{R}
\end{equation}
where $a>0$ is a constant such that
$\|{V}_{\delta} \|_2^2 =\int_{\RN}(I_{\alpha}*|{V}_{\delta}|^{\2})|{V}_{\delta}|^{\2}dx = S_{\alpha}^{\frac{N+\alpha}{\alpha}}$.
In the recent paper \cite{LiBao}, X. Li et al. also considered the minimizing problem (\ref{minimizingproblrm}). To settle the strict inequality (\ref{Theminiestimate}), they chose (\ref{extrmalLower}) as a test function and improved the existence result in \cite{Yao} to all $\mu>0$.

This problem seems to be more challenging when $ q \in( 2+\frac{4}{N}, 2^*)$ (the mass-supercritical region), since the functional $E_{\mu,c}$ is no longer bounded from below on $S(c)$.
In the present paper, we view $\mu>0$ and $c>0$ as two parameters, and throughout this paper,  for $p,q\in (2+\frac{4}{N},2^*)$, we always set
\begin{equation}\label{M}
\mathcal{M}_q(c,\mu):=\mu^{\frac{4}{N(q-2)-4}} c^{\frac{2\alpha}{N}+\frac{4(q-2)}{N(q-2)-4}}, \quad
\mathcal{N}_{p,q}(c,\mu):=[\mu c^{\frac{4}{N}}]^{\frac{4N(q-p)}{[N(q-2)-4]\cdot [N(p-2)-4]}}.
\end{equation}
Then it is easy to check that these mixed parameters satisfy the relation 
\begin{equation}
\frac{\mathcal{M}_p(c,\mu)}{\mathcal{M}_q(c,\mu)}=\mathcal{N}_{p,q}(c,\mu).
\end{equation}
Moreover, $\mathcal{M}_q(c,\mu)\to 0$ as $\mu\to 0$ (resp $c\to 0$) and  $\mathcal{M}_q(c,\mu)\to +\infty$ as $\mu\to \infty$ (resp $c\to \infty$).

In \cite{Yao}, S. Yao et al. also obtained
\begin{Prop}\label{Yao}{\cite{Yao}} Let $N \geq 2$, $\gamma>0$ and $p \in (2+\frac{4}{N},2^*)$. Then there exist $\mu_1 >0 $ such that for every $\mu>\mu_1$ and $c$ satisfying
\begin{equation}\label{yaocondition}
  \mathcal{M}_q(c,\mu) < \varpi:=\frac{2(N+\alpha)}{N}\cdot \frac{N(q-2)-4}{4q-2N(q-2)}\gamma^{-1} \|U\|_2^{\frac{4(q-2)}{N(q-2)-4}}S_{\alpha}^{\frac{N+\alpha}{N}},
\end{equation}
the problem $\left(\mathcal{A}_{\mu,c}\right)$ has a normalized ground state solution $(u,\lambda) \in S(c) \times \mathbb{R}^+$.
\end{Prop}
\noindent Here $U \in H_{rad}^1(\RN)$ denotes the unique solution \cite{Kwong} to
\begin{equation}\label{scalarequation}
  \left\{ \aligned  &-\Delta U + U =|U|^{q-2}U \quad \text{in\ \ }  \mathbb{R}^N,\\
    & U>0, \quad U(0)=\|U\|_{\infty}. \endaligned \right.
\end{equation}

Briefly speaking, following the idea in \cite{JeanjeanSIAM2019}, they constructed a sequence of radial function $\{ u_n \} \subset \mathcal{P}(\mu,c)$ and a $\{ \lambda_n \} \subset \mathbb{R}$ such that
$$\aligned & E_{\mu,c}(u_n) \rightarrow  \mathcal{E}(\mu,c):=\inf_{ v\in\mathcal{P}(\mu,c)}E_{\mu,c}(v), \quad E_{\mu,c}'(u_n) + \lambda_n u_n \rightarrow 0, \quad \lambda_n\rightarrow \bar{\lambda},\endaligned$$
where
$$\mathcal{P}(\mu,c):=\left\{ v\in S(c): \|\nabla v \|_2^2 =\mu \frac{N(q-2)}{2q}\|v \|_q^q\right\}.$$
Then they proved the sequence $\{ u_n \}$ converges strongly in $H^1(\RN)$ if $\mathcal{E}(\mu,c)>0$ and
\begin{equation}\label{Theestimate}
 \mathcal{E}(\mu,c) < \frac{\alpha}{2(N+\alpha)}\gamma^{-\frac{N}{\alpha}}(\bar{\lambda}S_{\alpha})^{\frac{N+\alpha}{\alpha}}-\frac{1}{2}\bar{\lambda}c^2.
\end{equation}
The assumptions (\ref{yaocondition}) and $\mu > \mu_1$ are crucial to their approach.
On the one hand, their compactness lemma relies heavily on $\E(\mu,c)>0$ which was guaranteed by
(\ref{yaocondition}). On the other hands, their argument to settle (\ref{Theestimate}) only works for $\mu > \mu_1$ with $\mu_1$ sufficiently large.

The purpose of this paper is two-fold.  
Firstly,  the following two problems will be addressed:
\begin{enumerate}
  \item [(1)] the existence of normalized solution under (\ref{yaocondition}) without the restriction $\mu > \mu_1$.
  \item [(2)] the existence of normalized solution for $\mathcal{M}_q(c,\mu)$ large enough.
\end{enumerate}

We point out that  (1) and (2) are largely open and the method employed  in \cite{Yao} seems to be 
fails to tackle (1) and (2).
To deal with those problems, we firstly introduce a novel compactness lemma and ascertain a new compactness threshold which can be strictly larger than $\frac{\alpha}{2(N+\alpha)}\gamma^{-\frac{N}{\alpha}}(\bar{\lambda}S_{\alpha})^{\frac{N+\alpha}{\alpha}}-\frac{1}{2}\bar{\lambda}c^2$.
To build the appropriate energy estimate, we present a novel approach to give a priori characterization of Lagrange multipliers.
It is noteworthy that the Lagrange multipliers were usually described only after the existence of normalized solutions was verified, so our method is significantly different from most of existing literatures. Furthermore, unlike the global branch approach developed by L.Jeanjean, J. Zhang and X. Zhong \cite{ZhongJMPA} nor the method proposed by L. Song and H. Hajaiej \cite{SongArxiv2022}, our method does not rely on the uniqueness or non-degeneracy of $U$, and the thresholds of $\mathcal{M}_q(c,\mu)$ will not be obtained by any limit process. In particular, if $g(u)=|u|^{q-2}u$, we prove that there exists two positive constants $\omega_1\le \omega_2$ such that Equ.$(\mathcal{A}_{\mu,c})$ admits a normalized ground state solution $(u_{\mu,c}, \lambda_{\mu,c})$
if $\mathcal{M}_q(c,\mu) \not\in [\omega_1,\omega_2]$.
The associated $(c,\mu)$ regions in which Equ.$(\mathcal{A}_{\mu,c})$ has a normalized ground state solution  and  a comparison with the previous results are depicted in the Figure 1 below.

Secondly, we study  the precise asymptotic profiles of normalized ground state solutions as the $(c, \mu)$-parameter varies.  In an excellent paper \cite{Wuyuanze}, Wei and Wu  investigate the existence and asymptotic behaviors of normalized solutions to a Schr\"odinger equation with combined powers nonlinearity, and obtain  a first  result on the precise asymptotic desscription of the normalized solutions based on ODE technique. However, in our setting,
due to  the presence of the nonlocal term, ODE technique does not work any more. 
In the present paper, after a suitable rescaling, we prove that
\begin{enumerate}
  \item [(1)]
  the solutions with $(c,\mu)\in {\bf (I)}$ converge to $U$ in $H^1(\RN)$ as  $\mathcal{M}_q(c,\mu)\rightarrow 0$, where  $U$ is the unique solution to (\ref{scalarequation});
   \item [(2)]
  the solutions with $(c,\mu)\in {\bf (II)}$ converge to $V_{\delta_0}$  in $H^1(\RN)$ as  $\mathcal{M}_q(c,\mu)\rightarrow +\infty$, where
 $V_{\delta_0}$ is a particular extremal function of $S_{\alpha}$ given in the next section.
 \end{enumerate}

$$\begin{tikzpicture}[samples=500]\label{Figue1}
 \node[above left]at(0,0){$0$};
\fill[gray!60](0,4)--plot[domain=0.25:5](\x,{1/\x})--(5,0.2)--(5,0)--(0,0);
\draw[dashed](0,4)--(0.25,4);
\draw[dashed](5,0.2)--(5,0);
 \draw[domain=1.25:5,dashed]plot(\x,{5/(\x)});
\fill[pattern=north east lines](0,4)--plot[domain=0.375:0.5](\x,{1.5/\x})--(0.75,2)--(0,2);
\fill[gray!60]plot[domain=1.25:5](\x,{5/\x})--(5,1)--(5,4)--(1.25,4);
 \draw[->](0,0)--(5.5,0) node[right]{$c$};
 \draw[->](0,0)--(0,4.5) node[left]{$\mu$};
  \fill(0,2)circle(1pt) node[left]{$\mu_1$};
  \fill(0,0)circle(1pt);
  \draw[domain=0.245:5,very thick]plot(\x,{1/(\x)});
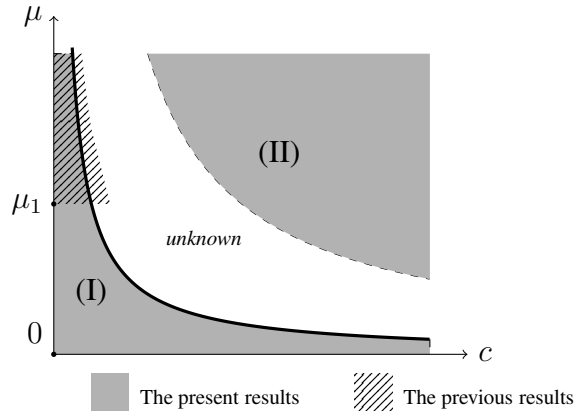
 \node[below] at (2.5,-1){Figure 1:  compare with the previous results};
 \fill[gray!60](0.5,-0.75)--(0.5,-0.25)--(1,-0.25)--(1,-0.75);
 \fill[pattern=north east lines](4,-0.75)--(4,-0.25)--(4.5,-0.25)--(4.5,-0.75);
  \node[right] at(4.5,-0.6){\tiny{The previous results}};
 \node[right] at (1,-0.6){\tiny{The present results}};
 \node[above] at (0.5,0.5){{(I)}};
  \node[below] at (3,3){{(II)}};
  \node[below] at (2,1.8){\tiny{\emph{unknown}}};\end{tikzpicture}$$

Finally, as a consequence of our main results, we present non-existence and multiplicity results of positive solutions to the problem $(\mathcal{B}_{\eta})$
\begin{equation*}
    \left\{ \aligned  &-\Delta u + u = \eta |u|^{q-2}u + (I_\alpha * |u|^{p})|u|^{p-2}u   & \text{in\ \ }  \mathbb{R}^N,\\
    & u \in H^1(\mathbb{R}^N), \endaligned \right. \eqno{(\mathcal{B}_{\eta})}
    \end{equation*}
where $N \geq 1$, $2+\frac{4}{N}\leq q<2^*$, $p=\frac{N+\alpha}{N}$ and  $\eta>0$.
A solution $v \in H^1(\RN)$ of $(\mathcal{B}_{\eta})$ can be characterized as a critical point of the action functional
\begin{equation}\label{Theactionfunctional}
  A_{\eta}(u):=\frac{1}{2}\|\nabla u\|_2^2 + \frac{1}{2}\| u \|_2^2 - \frac{\eta}{q}\| u\|_q^q - \frac{N}{2(N+\alpha)}\II,
\end{equation}
and we say a nontrivial solution $v$ is least action solution if it achieves the infimum of $A_{\eta}$ among all the nontrivial solutions, namely,
$$A_{\eta}(v)= a(\eta):=\inf\left\{  A_{\eta}(w): w\in H^1(\RN)\setminus\{0\}:\    A_{\eta}'(w)=0    \right\}. $$
The existing literatures commonly focused on the existence of least action solutions, and the following result was established in \cite{LiMa,Tangxianhua1}.
\begin{Prop}\label{resultCCM}
Let $N \geq 1$, $\alpha \in (0,N)$ and $\eta>0$. Then there exists a constant $\eta_0>0$ such that $(B_{\eta})$ admits a positive least action solution $v_{\eta} \in H^1(\RN)$
which is radially symmetric and radially nonincreasing if one of the following conditions holds
\begin{enumerate}
  \item [(1)] $q \in (2,2+\frac{4}{N})$;
  \item [(2)] $q \in [2+\frac{4}{N}, \frac{2N}{N-2})$ and $\eta > \eta_0$.
\end{enumerate}
\end{Prop}
In this paper, we find a threshold $\eta_1 \leq \eta_0$ such that
$(\mathcal{B}_{\eta})$ admits no least action solution for $\eta<\eta_1$, and admits a least action solution for $\eta \geq \eta_1$.  A further surprising result is that we find a constant $\eta_2 \geq \eta_1$ such that $(\mathcal{B}_{\eta})$ admits two positive radial solutions for $\eta>\eta_2$ in the case $2+\frac{4}{N}<q<2^*$. This seems to be the first result concerning the non-existence and multiplicity of positive solutions for Choquard type equations involving the lower critical exponent.

We remark that the multiplicity of positive solutions for elliptic equations in the whole space is delicate and difficult. By using  Lyapunov-Schmidt type arguments,  D\'avila, del Pino and Guerra \cite{Davila-1} constructed three positive solutions of the equation
\begin{equation}\label{e117}
-\Delta v+v=v^{p-1} + \eta v^{q-1}\quad  {\rm in} \ \mathbb R^3
\end{equation}
for  sufficiently large $\eta>0$, $2<q<4$  and slightly subcritical $p < 6$.  Wei and Wu \cite{Wei-2} obtained two positive solutions for the equation \eqref{e117} with $p=6$, $2<q<4$ and sufficiently large $\eta>0$. Very recently, when $N=3$, the second author \cite{Ma} study the Choquard type equation $(\mathcal{B}_{\eta})$ involving the upper exponent $p=\frac{N+\alpha}{N-2}$ and obtain two positive solutions for $2<q<4$ and sufficiently large $\eta>0$.
Our method is based on a novel idea and quite different from the previous related results \cite{Davila-1,Wei-2,Ma} in nature.

\noindent \textbf{Organization of the paper.} In section 2, we present the main results. In Section 3, we give some preliminary results.  In Section 4, a new
minimax theorem and a new compactness lemma are presented. In Section 5, by introducing some new estimates, we show the existence of ground state solutions and prove Theorems \ref{TheSecondResult} and \ref{theresult2}. The Section 6 is devoted to study the
asymptotic profiles of the solutions where Theorems \ref{ap1}--\ref{ap2} are proved. In Section 7, we consider the power type nonlinearity and prove Theorems  \ref{Homoresult} --\ref{Mutli}.

\noindent \textbf{Notations.} From now on in this paper, otherwise mentioned, we use the following notations:
\begin{enumerate}
  \item [$\bullet$] $L^s(\mathbb{R}^N)$ with $ s \in [1,\infty)$ is the Lebesgue space with the norm $\|u \|_s=(\int_{\mathbb{R}^N}|u|^s dx )^{1/s}$.
  \item [$\bullet$] $H^{1}(\mathbb{R}^{N})$ is the usual Sobolev space with the norm
$\| u \|_{H^1}=\left( \int_{\mathbb{R}^N}|\nabla u|^2+|u|^2 dx\right)^{1/2}$.
  \item [$\bullet$] $H_{rad}^1(\mathbb{R}^N)$ denotes $\{ u \in H^1(\mathbb{R}^N):u\ \text{is\ radial\ symmetric}\}$.
  \item [$\bullet$] $C_1$, $C_2, \dots$ denote positive constants possibly different in different places.
  \item [$\bullet$] For any $u \in H^1(\RN)$ and $t>0$, $u_t(x)=t^{\frac{N}{2}}u(tx)$, $\forall x \in \RN$.
  \item [$\bullet$] $u_n \rightharpoonup u$ in $H^1(\RN)$ denotes $\{u_n\}$ converges weakly to $u$ in $H^1(\RN)$.
\end{enumerate}
For any $s>0$ and two nonnegative functions $f(s)$ and $h(s)$, we write
\begin{enumerate}
  \item [(1)] $f(s) \lesssim h(s)$ or $h(s) \gtrsim f(s)$ if there exists a positive constant $C$ independent of $s$ such that $f(s) \leq C h(s)$.
  \item [(2)] $f(s) \thicksim h(s)$ if $f(s) \lesssim h(s)$ and $f(s) \gtrsim h(s)$.
\end{enumerate}

\section{Main Results}

To state our main results, we make the following  assumptions
\begin{enumerate}
  \item [(G1)] $g: \mathbb{R}\rightarrow \mathbb{R}$ is continuous and odd.
  \item [(G2)] There exist $2+\frac{4}{N} < p\leq q < 2^*(+\infty,N=1,2)$ such that $ 0< pG(t) \leq g(t)t \leq qG(t)$.
  \item [(G3)] Let $\bar{G}(t):=g(t)t-2G(t)$, the function $t \mapsto \bar{G}(t)/|t|^{2+\frac{4}{N}}$ is strictly decreasing on $(-\infty,0)$ and strictly increasing on $(0,+\infty)$.
  \item [(G4)] There exists a positive constant $A_1$ such that
                   $\lim_{t \rightarrow +\infty} \frac{g(t)}{|t|^{q-1}}=A_1.  $
  \item [(G5)] There exists a positive constant $A_2$ such that
                   $\lim_{t \rightarrow 0} \frac{g(t)}{|t|^{p-1}}=A_2.  $
\end{enumerate}
\begin{Rem} Particularly, the function $g(t)=\Sigma_{i=1}^{k}a_i |t|^{q_i-2}t$ with $a_i>0$ and $q_i \in (2+\frac{4}{N},2^*)$ satisfies the assumptions (G1)--(G5).
Following \cite[(2.2)]{Jeanjean1997} we obtain for all $s \in \mathbb{R}$,
\begin{equation}\label{barGGg}
  \left\{ \aligned  &  \frac{1}{q-2} \bar{G}(s) \leq G(s) \leq \frac{1}{p-2} \bar{G}(s)    \\
                     &  \frac{q}{q-2} \bar{G}(s) \leq g(s)s \leq \frac{p}{p-2} \bar{G}(s).     \endaligned \right.
\end{equation}
\end{Rem}
In Theorems 2.1-2.4, we always assume that  $p<q$, our first two results concern the existence of normalized solutions. 
\begin{Theorem}\label{TheSecondResult}
Let $N\geq 1$, $\gamma>0$ and $g$ satisfies (G1)--(G4), then there exists a positive constants $\delta$ independent of $\mu,c$ such that for every $\mu$, $c$ satisfying
$$\max \left\{\mathcal{M}_q(c,\mu), \mathcal{N}_{p,q}(c,\mu)\right\} < \delta, $$
Equ. (\ref{ThePro}) possesses a normalized ground state solution $({u}_{\mu,c},{\lambda}_{\mu,c})$.
Moreover, ${u}_{\mu,c}$ is a positive radial, nonincreasing function of $|x|$, 
\begin{equation}\label{TheCharOfLag}
{\lambda}_{\mu,c}\sim  \mu^{-\frac{4}{N(q-2)-4}} c^{-\frac{4(q-2)}{N(q-2)-4}}
\end{equation}
and
\begin{equation}\label{TheCharOfEnergy}
 E_{\mu,c}(u_{\mu,c}) + \frac{N\gamma}{2(N+\alpha)}c^{\frac{2(N+\alpha)}{N}}S_{\alpha}^{-\frac{N+\alpha}{N}} \sim \mu^{-\frac{4}{N(q-2)-4}} c^{-\frac{4N-2q(N-2)}{N(q-2)-4}}.
\end{equation}
\end{Theorem}


\begin{Theorem}\label{theresult2}
Let $N\geq 1$, $\gamma>0$ and $g$ satisfies (G1)--(G3) and (G5), then there exists a positive constant $\tau$ independent of $\mu,c$ such that
for every $\mu$, $c$ satisfying
$$\min\left\{\mathcal{M}_q(c,\mu), \mathcal{N}_{p,q}(c,\mu)\right\}>\tau, $$
Equ. (\ref{ThePro}) possesses a normalized ground state solution $({u}_{\mu,c},{\lambda}_{\mu,c})$.
Moreover, ${u}_{\mu,c}$ is a positive, radial, nonincreasing function of $|x|$,
\begin{equation}\label{TheCharLag2}
  \Big| \lambda_{\mu,c}c^2 -\gamma S_{\alpha}^{-\frac{N+\alpha}{N}}c^{\frac{2(N+\alpha)}{N}}\Big| \lesssim \mu^{-\frac{4N-2p(N-2)}{N(p-2)-4}} c^{-\frac{4(p-2)}{N(p-2)-4}}
\end{equation}
and
\begin{equation}\label{TheCharOfEnergy2}
  E_{\mu,c}(u_{\mu,c}) + \frac{N\gamma}{2(N+\alpha)}c^{\frac{2(N+\alpha)}{N}}S_{\alpha}^{-\frac{N+\alpha}{N}} \sim \mu^{-\frac{4}{N(p-2)-4}}  c^{-\frac{4N-2p(N-2)}{N(p-2)-4}}.
\end{equation}
\end{Theorem}

The following two results are  devoted to 
the asymptotic profiles of $(u_{\mu,c}, \lambda_{\mu,c})$ obtained above.
\begin{Theorem}\label{ap1}
Let $A_1=1$ and $({u}_{\mu,c}, {\lambda}_{\mu,c})$ be a family of normalized ground state solutions to \eqref{ThePro} obtained in Theorem \ref{TheSecondResult}, then
$$ \mu^{\frac{N}{N(q-2)-4}}c^{\frac{4}{N(q-2)-4}}u_{\mu,c}\left(\mu^{\frac{2}{N(q-2)-4}}c^{\frac{2(q-2)}{N(q-2)-4}}x\right)  \rightarrow  \|U\|_2^{\frac{4}{N(q-2)-4}}U\left(\|U\|_2^{\frac{2(q-2)}{N(q-2)-4}}x\right)$$
in $H^1(\RN)$
and
$$\mu^{\frac{4}{N(q-2)-4}} c^{\frac{4(q-2)}{N(q-2)-4}}\lambda_{\mu,c} \rightarrow \|U\|_{2}^{\frac{4(q-2)}{N(q-2)-4}} $$
as $\max\left\{\mathcal{M}_q(c,\mu), \mathcal{N}_{p,q}(c,\mu)\right\}\rightarrow 0$.
\end{Theorem}

\begin{Theorem}\label{ap2}
  Let $\gamma=1$, $A_2=1$, and $({u}_{\mu,c}, {\lambda}_{\mu,c})$ be a family of normalized ground state solutions to $(\ref{ThePro})$ obtained in Theorem \ref{theresult2}, then there exists a family $\{ \zeta_{\mu,c} \} \subset \mathbb{R}$ such that
  $$\zeta_{\mu,c}^{\frac{N}{2}}u_{\mu,c}(\zeta_{\mu,c} x) \rightarrow V_{\delta_0}  $$
  in $H^1(\RN)$
  as $\min\left\{ \mathcal{M}_p(c,\mu), \mathcal{N}_{p,q}(c,\mu)\right\}\rightarrow +\infty$, where
$$ \zeta_{\mu,c}^2 \sim \mathcal{M}_p(c,\mu) $$
and
$$  {V}_{\delta_0}(x)=a \left( \frac{\delta_0}{\delta_0^2 + |x|^{2}} \right)^{\frac{N}{2}}$$
is an extremal function  of $S_\alpha$ such that $\|V_{\delta_0}\|_2^2 =S_{\alpha}^{\frac{N+\alpha}{\alpha}}$ and $\|\nabla V_{\delta_0} \|_2^2=\frac{N(p-2)}{2p}\| V_{\delta_0}\|_p^p$.
\end{Theorem}

In the following, we assume that $g(u)=|u|^{q-2}u$, then the following result may be obtained. 
\begin{Theorem}\label{Homoresult} Let $N \geq 1$, $\gamma > 0$, $\mu>0$, $c>0$, and $2+ \frac{4}{N} < q <2^* $. Then the followings are true:\\
\noindent(1)\ (Existence results). There exist two positive constants $\omega_1$ and $\omega_2$ independent of $\mu,c$ such that for every $\mu$, $c$ satisfying
\begin{equation*}
  \mathcal{M}_q(c,\mu) \in (0,\omega_1) \cup (\omega_2,+\infty),
\end{equation*}
Equ.$(\mathcal{A}_{\mu,c})$ admits a normalized ground state solution $(u_{\mu,c}, \lambda_{\mu,c})$.
Moreover, ${u}_{\mu,c}$ is a positive, radial, nonincreasing function of $|x|$.

\noindent(2)\ (Asymptotic profiles). Let $\mathcal{M}_q(c,\mu)\rightarrow 0$, then
$$ \mu^{\frac{N}{N(q-2)-4}}c^{\frac{4}{N(q-2)-4}}u_{\mu,c}\left(\mu^{\frac{2}{N(q-2)-4}}c^{\frac{2(q-2)}{N(q-2)-4}}x\right)  \rightarrow  \|U\|_2^{\frac{4}{N(q-2)-4}}U\left(\|U\|_2^{\frac{2(q-2)}{N(q-2)-4}}x\right)$$
in $H^1(\RN)$
and
$$\mu^{\frac{4}{N(q-2)-4}} c^{\frac{4(q-2)}{N(q-2)-4}}\lambda_{\mu,c} \rightarrow \|U\|_{2}^{\frac{4(q-2)}{N(q-2)-4}}. $$
Let $\gamma=1$ and $\mathcal{M}_q(c,\mu) \rightarrow +\infty$, then
there exists a family $\{ \zeta_{\mu,c} \} \subset \mathbb{R}$ such that
  $$\zeta_{\mu,c}^{\frac{N}{2}}u_{\mu,c}(\zeta_{\mu,c} x) \rightarrow V_{\delta_0}  $$
  in $H^1(\RN)$,
 where
$$ \zeta_{\mu,c}^2 \sim \mathcal{M}_q(c,\mu) $$
and $V_{\delta_0}$ is an extremal function of $S_\alpha$ such that $\|V_{\delta_0}\|_2^2 =S_{\alpha}^{\frac{N+\alpha}{\alpha}}$ and $\|\nabla V_{\delta_0} \|_2^2=\frac{N(q-2)}{2q}\| V_{\delta_0}\|_q^q$.
\end{Theorem}
\begin{Rem} In particular, Theorem \ref{Homoresult} implies that  for any $\mu>0$ fixed, there exist  $0<c_1\le c_2<\infty$, such that  Equ. $(\mathcal{A}_{\mu,c})$ admits a normalized ground state solution if $c\notin [c_1,c_2]$. Thus we improve sharply  the  previous result \cite{Yao}.


\end{Rem}

Inspired by the above facts, we shall explore the non-existence and multiplicity of positive radial solutions to the following problem
\begin{equation*}
    \left\{ \aligned  &-\Delta u + u = \eta |u|^{q-2}u + (I_\alpha * |u|^{\frac{N+\alpha}{N}})|u|^{\frac{N+\alpha}{N}-2}u   & \text{in\ \ }  \mathbb{R}^N,\\
    & u \in H^1(\mathbb{R}^N), \endaligned \right. \eqno{(\mathcal{B}_{\eta})}
    \end{equation*}
where $N\ge 1$ and $2+\frac{4}{N} \leq q<2^*$, then we have the following result.

 \begin{Theorem}\label{Mutli}
Let $N \geq 1$, $2+\frac{4}{N} \leq q<2^*$, then there exist $\eta_1 \leq \eta_2$ such that
\begin{enumerate}
  \item [(1)] $(\mathcal{B}_{\eta})$ admits no positive radial least action solution for $\eta < \eta_1 $;
  \item [(2)] $(\mathcal{B}_{\eta})$ admits a positive radial least action solution $v_{\eta}$ for $\eta \geq \eta_1$ in the case of $2+\frac{4}{N} < q<2^*$, and for $\eta > \eta_1$ in the case of $q=2+\frac{4}{N}$.
  Moreover, $A_{\eta}(v_{\eta}) \rightarrow 0$ as $\eta \rightarrow +\infty$;
  \item [(3)]  $(\mathcal{B}_{\eta})$ admits a positive radial solutions $u_{\eta}$ differs from $v_{\eta}$ for $\eta > \eta_2$ in the case of $2+\frac{4}{N} < q<2^*$. Moreover, $A_{\eta}(u_{\eta}) \rightarrow \frac{\alpha}{2(N+\alpha)}S_{\alpha}^{\frac{N+\alpha}{\alpha}}$ as $\eta \rightarrow +\infty$.
\end{enumerate}
\end{Theorem}
\begin{Rem}
Proposition \ref{resultCCM} (1), Theorem \ref{Mutli} (1) and (2) almost completely solves the existence of least action to $(\mathcal{B}_{\eta})$ expect for $q=2+\frac{4}{N}$ and $\eta=\eta_1$. Theorem \ref{Mutli} (3) implies the multiplicity of positive radial solutions to $(\mathcal{B}_{\eta})$ for $\eta$ sufficiently large.
\end{Rem}

Now let us sketch our proof and highlight some of the difficulties encountered. We will use a natural constraint, namely the Pohoz\v{a}ev manifold
\begin{equation}\label{TheManifold}
  \mathcal{P}(\mu,c)=\left\{ u\in S(c): P_{\mu,c}(u): =\|\nabla u \|_2^2 - \mu \frac{N}{2}\int_{\RN}\bar{G}(u)dx =0  \right\},
\end{equation}
and search the minimizer of \
\begin{equation}\label{TheMinimizingProblem}
  \mathcal{E}(\mu,c):=\inf_{v \in \mathcal{P}(\mu,c)}E_{\mu,c}(v).
\end{equation}
Since $\mathcal{P}(\mu,c)$ contains all nontrivial solutions to (\ref{ThePro}), the minimizer $u_{\mu,c}$ of $\E(\mu,c)$ is clearly a normalized ground state solution.

Motivated by \cite{JeanjeanSIAM2019} and \cite{Yao}, we construct a bounded Palais-Smale sequence $\{ u_n \}$ at the level $\E(\mu,c)$.
It is worth to point out that the existence result is open for the case $N=1$ while \cite{Yao} only consider $N \geq 2$ using the compactly embedding $H^1_{rad}(\RN) \hookrightarrow L^r_{rad}(\RN)$ for $ r \in (2, 2^*)$.
In order to deal with the case $N=1$,
we present a slightly different minimax argument involved with the Schwartz rearrangement.
This general minimax theorem establish a bounded sequence $\{\tilde{u}_n \}\subset \mathcal{P}(\mu,c)$ with $\tilde{u}_n(x)$ is radial and nonincreasing function of $|x|$ for every $n \geq 0$, such that $\| u_n -\tilde{u}_n \|_{H^1(\RN)} \rightarrow 0 $ as $n \rightarrow \infty$. Therefore, it is sufficient to consider the convergence of $\{u_{n}\}$ since $\{ \tilde{u}_n \}$ has a strongly convergent subsequence in $L^r(\RN)$ for $N \geq 1$ and $r \in (2, 2^*)$.
This idea is credited to \cite[Lemma 4.2, Remark 5.2]{SoaveJDE} and \cite[Proposition 1.7.1]{TheSemilinearSchrodinger}.

Adopting the Lagrange multipliers rule, we obtain a sequence $\{\lambda_n \} \subset \mathbb{R}$, then we may assume
$u_n \rightharpoonup {u}_{\mu,c}$ weakly in $H^1(\RR)$ and $\lambda_n \rightarrow {\lambda}_{\mu,c}$ in $\mathbb{R}$.
To establish the compactness one needs to obtain a rigorous upper bound estimates for the energy levels
while excluding the ``bubbling phenomenon" due to the nonlocal term $\II$.
Our strategy is a new compactness lemma (see Lemma \ref{compactnesslemma}) which is one of the major novelties.
Let $w_n=u_n-u_{\mu,c}$, then $\{u_n\}$ will split into two parts $u_{\mu,c}$ and $\{w_n\}$. By a standard argument based on the Brezis-Lieb type lemma and (\ref{TheLowerConstant}) one has
$$\aligned & E_{\mu,c}(u_n) + \frac{1}{2} \lambda_n \| u_n\|_2^2 +o(1)\\
= & E_{\mu,c}(u_{\mu,c})+ \frac{\lambda_{\mu,c}}{2}\|u_{\mu,c} \|_2^2 + \frac{\lambda_n}{2}\|w_n\|_2^2 - \frac{N\gamma}{2(N+\alpha)}\int_{\RN}(I_{\alpha}*|w_n|^{\2})|w_n|^{\2}dx\\
\geq & E_{\mu,c}(u_{\mu,c})+ \frac{\lambda_{\mu,c}}{2}\|u_{\mu,c} \|_2^2 +\frac{\alpha}{2(N+\alpha)}\gamma^{-\frac{N}{\alpha}}(\lambda_{\mu,c} S_{\alpha})^{\frac{N+\alpha}{\alpha}}. \endaligned $$
Notice that if $u_{\mu,c} \not\equiv 0$ and $\lambda_{\mu,c}>0$, then $u_{\mu,c}$ is a nontrivial solution to
\begin{equation}\label{Thefixedfrequencygeneralintroduction}
   -\Delta {u} + {\lambda}_{\mu,c} {u} =\mu g(u) + \gamma (I_{\alpha}*|u|^{\2})| u|^{\2-2} u \quad \text{\ in\ }\RN,
\end{equation}
and hence, we have $w_n \rightarrow 0$  and $u_n \rightarrow u_{\mu,c}$ in $H^1(\RN)$ if
\begin{equation}\label{strictinequalityintroduction}
 \mathcal{E}(\mu,c)+\frac{1}{2}\lambda_{\mu,c}c^2 < a(\mu,c)+\frac{\alpha}{2(N+\alpha)}\gamma^{-\frac{N}{\alpha}}(\lambda_{\mu,c} S_{\alpha})^{\frac{N+\alpha}{\alpha}}
\end{equation}
where $a(\mu,c)$ is the least action level of problem (\ref{Thefixedfrequencygeneralintroduction}).
This new compactness threshold can be strictly large than that used in \cite{Yao} as $a(\mu,c) > 0$ in general.

The primary difficulty in establishing (\ref{strictinequalityintroduction}) lies in providing an a priori characterization of $ \lambda_{\mu,c}$.
Most of the existing literature describe the multipliers only after the existence of solutions has been verified. Therefore, novel method must be developed to tackle the problem, particularly since the ODE approach in \cite{Wuyuanze} fails due to the presence of nonlocal term.
 By employing the Pohoz\v{a}ev identity,
the relationship between $\lambda_{\mu,c}$ and $\mathcal{E}(\mu,c)$ is concluded. Furthermore, by selecting an appropriate test function, we obtain a robust estimate of $\mathcal{E}(\mu,c)$ and  a precise characterization of $\lambda_{\mu,c}$ as a consequence. Through some detailed calculations, we find two positive constants $\omega_1\le \omega_2$ such that (\ref{strictinequalityintroduction}) hold if $\mathcal{M}_q(c,\mu) \not\in [\omega_1,\omega_2]$ and this complete the proof.
Base on those detailed estimates, we also derive the sharp quantitative properties of $(u_{\mu,c},\lambda_{\mu,c})$ as the parameter $(c,\mu)$ varies.

Then let $\{(u_{\mu}, \lambda_{\mu})\}$ denotes the family of the normalized ground state solutions to $\left(\mathcal{A}_{\mu,c}\right)$ obtained in Theorem \ref{Homoresult} for $c>0$ fixed, $\gamma=1$ and $\mu \rightarrow +\infty$,
by introducing rescaling by
$$ \bar{u}_{\mu}:= \lambda_{\mu}^{-\frac{N(2+\alpha)}{4\alpha}}u_{\mu}(\lambda_{\mu}^{-\frac{1}{2}}x),      $$
we know that $\bar{u}_{\mu}$ solves
\begin{equation*}
    \left\{ \aligned  &-\Delta u + u = \eta_{\mu} |u|^{q-2}u + (I_\alpha * |u|^{\frac{N+\alpha}{N}})|u|^{\frac{N+\alpha}{N}-2}u   & \text{in\ \ }  \mathbb{R}^N,\\
    & u \in H^1(\mathbb{R}^N) \endaligned \right.
    \end{equation*}
where $\eta_{\mu} \rightarrow +\infty$ as $\mu \rightarrow +\infty$. When $\mu$ is sufficiently large, it is deduced that $\bar{u}_{\mu}$ is distinct from the least action solution, as evidenced by these quantitative properties of $\{(u_{\mu}, \lambda_{\mu})\}$.

\section{Preliminaries}

The Gagliardo-Nirenberg inequality can be seen in \cite{Weinstein}
\begin{Lem}\label{GNS}
Let $N \geq 1$, $2< s <2^*$ and $\gamma_s =\frac{N(s-2)}{2s} $, then the following sharp inequality
$$\| u \|_s \leq C_{N,s}\| u\|_2^{1-\gamma_s} \| \nabla u\|_2^{\gamma_s}    $$
holds for any $u \in H^1(\RN)$, where the constant $C_{N,s}=\left( \frac{s}{  2\| Q_s \|_2^{s-2}  }  \right)^{\frac{1}{s}}$
and $Q_s$ is the unique positive solution of
\begin{equation}\label{equ:Q sovles equ}
  -\frac{N(s-2)}{4} \Delta Q_s + \frac{2N-s(N-2)}{4} Q_s=|Q_s|^{s-2}Q_s, \ \ x\in \mathbb{R}^N.
\end{equation}
\end{Lem}
\begin{Rem}
$$C_{N,q}^q =\frac{2q}{2q-N(q-2)}\left(  \frac{2q-N(q-2)}{N(q-2)}\right)^{\frac{N(q-2)}{4}} \| U \|_2^{2-q}.  $$
\end{Rem}

The Hardy-Littlewood-Sobolev inequality can be seen in \cite[Theorem 4.3]{Lieb2001}

\begin{Lem}
Let $N \geq 1$, $s_1$, $s_2 \geq 1$, and $0 < \beta < N$ with $\frac{1}{s_1}+\frac{(N-\beta)}{N}+\frac{1}{s_2}=2$. Let $u \in L^{s_1}(\RN)$ and $v \in L^{s_2}(\RN)$. Then there exists a sharp constant $C(N, \beta, s_1)$, independent of $u$ and $v$, such that
$$ \Biggl| \int_{\RN}\int_{\RN}   \frac{u(x)v(y)}{|x-y|^{N-\beta}}dxdy    \Biggr|  \leq C(N, \beta, s_1)\|u \|_{s_1} \|v \|_{s_2}.   $$
\end{Lem}
\begin{Lem}\label{HLWSE}
Let $N \geq 1$, $\alpha \in (0,N)$, and $S_\alpha$ be defined in (\ref{TheLowerConstant}). Then $S_{\alpha}>0$ and the minimizer of $S_{\alpha}$ is
$${V}_{\delta}(x)=a \left( \frac{\delta}{\delta^2 + |x-y|^{2}}   \right)^{\frac{N}{2}}, \quad \delta>0,\quad y \in\mathbb{R}$$
where $a>0$ is a constant such that
$$\|{V}_{\delta} \|_2^2 =\int_{\RN}(I_{\alpha}*|{V}_{\delta}|^{\2})|{V}_{\delta}|^{\2}dx = S_{\alpha}^{\frac{N+\alpha}{\alpha}}. $$
\end{Lem}


The following Br\'{e}zis-Lieb-type lemma can be found in \cite[Lemma 2.4]{MorozJFA}.
\begin{Lem}
Let $N \geq 1$, $\alpha \in (0,N)$ and $\{ u_n \}$ be a bounded sequence in $L^2(\RN)$. If $u_n \rightarrow u$ a.e. in $\RN$, then
$$\aligned \IIN  =& {\int_{\RN}(I_{\alpha}*|u_n-u|^{\2})|u_n-u|^{\2}dx} \\
                  &+  \II + o_n(1). \endaligned    $$
\end{Lem}

The following compactness lemma can be found in \cite[Proposition 1.7.1]{TheSemilinearSchrodinger}.
\begin{Lem}\label{TheCompactnessTool}
Let $\{ u_n\} \subset H^1(\RN)$ be a bounded sequence of spherically symmetric functions. If $N \geq 2$ or if $u_n(x)$ is a nonincreasing function of $|x|$ for every $n \geq 0$, then, up to a subsequence,
there exist $u\in H^1(\RN)$ such that $u_n \rightarrow u$ in $L^r(\RN)$ for every $2<r<2^*$($2<r\leq \infty $ if $N=1,2$).
\end{Lem}

 For the homogeneous nonlinearity $g(u)= |u|^{q-2}u$, we define $e(\mu,c):=\inf_{u \in \mathcal{P}(\mu,c)} I(u)$, and $I(u):= \frac{1}{2}\|\nabla u\|_2^2 - \frac{1}{q}\| u\|_q^q $

\begin{Lem}\label{Them(mu,c)}
Let $N \geq 1$, $2+\frac{4}{N} < q < 2^*$. Then for every $\mu$, $c>0$, $e(\mu,c)$ admits a unique (up to translations) minimizer $(Z_{\mu,c}, {\Lambda}_{\mu,c})$, and
\begin{equation}\label{Thefo}
    \aligned &Z_{\mu,c}(x):= \mu^{-\frac{1}{q-2}}\rho^{-\frac{2}{q-2}}U(x/\rho), &{\Lambda}_{\mu,c}= \rho^{-2} \endaligned
\end{equation}
where $\rho = \mu^{\frac{2}{N(q-2)-4}} c^{-\frac{2(q-2)}{N(q-2)-4}} \|U\|_2^{\frac{2(q-2)}{N(q-2)-4}}$.
\end{Lem}
\Proof Let $U \in H^1_{rad}(\RN)$ be the unique solution to (see\cite{Kwong})
$$ \left\{ \aligned  &-\Delta U + U =|U|^{q-2}U \quad \text{in\ \ }  \mathbb{R}^N,\\
    & U>0, \quad U(0)=\|U\|_{\infty}. \endaligned \right.$$
Then for $c>0$ and $\mu>0$ one obtains a solution $(Z_{\mu,c}, {\Lambda}_{\mu,c}) \in H^1_{rad}(\RN) \times \mathbb{R}^+$ of
\begin{equation}\label{Minimizerequation}
  \left\{ \aligned  &-\Delta Z_{\mu,c} + {\Lambda}_{\mu,c} Z_{\mu,c} =\mu |Z_{\mu,c}|^{q-2}Z_{\mu,c} \quad \text{in\ \ }  \mathbb{R}^N,\\
    & Z_{\mu,c}>0, \quad \|Z_{\mu,c} \|_2^2=c^2 \endaligned \right.
\end{equation}
by scaling:
\begin{equation*}
  \aligned &Z_{\mu,c}(x):= \mu^{-\frac{1}{q-2}}\rho^{-\frac{2}{q-2}}U(x/\rho), &{\Lambda}_{\mu,c}= \rho^{-2} \endaligned
\end{equation*}
where $\rho $ is determined by
$$ c^2 =\|Z_{\mu,c}\|_2^2=\rho^{N-\frac{4}{q-2}}\mu^{-\frac{2}{q-2}}\rho_0^2, \quad \rho_0=\|U \|_2.    $$
This completes the proof. \qed

\begin{Lem}\label{TheMaxFiber}
Let $N \geq 1$, $\mu>0$, $c>0$ and $g$ satisfy (G1)--(G3). The followings are true:
\begin{enumerate}
  \item [(1)] For any $u \in S(c)$,
there exists a unique $t(u)>0$ such that $u_{t(u)} \in \mathcal{P}(\mu,c)$ and
$E_{\mu,c}(u_{t(u)})=\max_{t>0}E_{\mu,c}(u_t)$;
  \item [(2)] $\E(\mu,c):= \inf_{u \in \mathcal{P}(\mu,c)}E_{\mu,c}(u)= \inf_{u \in S(c)}\max_{t>0}E_{\mu,c}(u_t) $;
  \item [(3)] The mapping $u \mapsto t(u)$ is continuous in $ u \in H^1(\RN)\setminus \{ 0\}$;
  \item [(4)] $E_{\mu,c}$ is coercive and bounded from below on $\mathcal{P}(\mu,c)$;
  \item [(5)] If suppose further $(G4)$ hold, then $$ \E(\mu,c) + \frac{N\gamma}{2(N+\alpha)}c^{\frac{2(N+\alpha)}{N}}S_{\alpha}^{-\frac{N+\alpha}{N}} \lesssim  \mu^{-\frac{4}{N(q-2)-4}} c^{-\frac{4N-2q(N-2)}{N(q-2)-4}}.$$
\item [(6)] If suppose further $(G5)$ hold, then $$\E(\mu,c) + \frac{N\gamma}{2(N+\alpha)}c^{\frac{2(N+\alpha)}{N}}S_{\alpha}^{-\frac{N+\alpha}{N}} \lesssim  \mu^{-\frac{4}{N(p-2)-4}} c^{-\frac{4N-2p(N-2)}{N(p-2)-4}}.$$
\end{enumerate}
\end{Lem}
\Proof (1) (2) (3) The proofs are similar to \cite[Lemma 2.4]{JeanjeanCV2020} and \cite[Lemma 2.10]{Jeanjean1997} and omitted.

(4)\ It follows from (G2)(G3) that
\begin{equation}\label{Edec}
  E_{\mu,c}(u) \geq \left(\frac{N(p-2)-4}{2N(p-2)} \right)  \|\nabla u\|_2^2 -  \frac{\gamma N}{2(N+\alpha)}S_{\alpha}^{-\frac{N+\alpha}{N}}c^{\frac{2(N+\alpha)}{N}}
\end{equation}
for any $u \in \mathcal{P}(\mu,c)$
and hence, $E_{\mu,c}$ is coercive and bounded from below on $\mathcal{P}(\mu,c)$.

(5) It follows from (G2) and (G4) that there exists $C_1>0$ such that $G(t) \geq  \frac{C_1}{q}|t|^{q}$ for any $t \in \mathbb{R}$. Then one has
$$E_{\mu,c}(u) \leq I_{\mu,c}(u):=\frac{1}{2}\|\nabla u\|_2^2 -\frac{\mu}{q}C_1 \|u \|_q^q -\II  $$
for any $u\in H^1(\RN)$.
Let $v_1(x) = c{S_{\alpha}^{-\frac{N+\alpha}{\alpha}}}{V}_{1}(x)$ and ${V}_{1}(x)$ was defined in Lemma \ref{HLWSE}, a direct calculation shows that
\begin{equation}\label{testbyV1}
  \|\nabla v_1\|_2^2=c^2 S_{\alpha}^{-\frac{2(N+\alpha)}{\alpha}} \|\nabla {V}_{1} \|_2^2,\ \|v_1\|_2^2=c^2, \ \|v_1 \|_q^q =c^q S_{\alpha}^{-\frac{q(N+\alpha)}{\alpha}} \|{V}_{1}\|_q^q,
\end{equation}
and
\begin{equation}\label{testbyV2}
\int_{\RN}(I_{\alpha}*|v_1|^{\2})|v_1|^{\2}dx = c^{\frac{2(N+\alpha)}{N}}S_{\alpha}^{-\frac{N+\alpha}{N}}.
\end{equation}
By the definition of $\E(\mu,c)$ and (\ref{testbyV1}),(\ref{testbyV2}), we get
\begin{equation}\label{UpperNewcalcu}
 \aligned  \E(\mu,c) \leq \max_{t>0}I_{\mu,c}((v_1)_t) =& C \mu^{-\frac{4}{N(q-2)-4}} c^{\frac{2N(q-2)-4q}{N(q-2)-4}} \\ &-\frac{N\gamma}{2(N+\alpha)}c^{\frac{2(N+\alpha)}{N}}S_{\alpha}^{-\frac{N+\alpha}{N}} \endaligned
\end{equation}
for any $\mu>0$ and $c>0$. \qed

(6)\ It follows from (G2) and (G5) that there exists $C_2>0$ such that $G(t) \geq  \frac{C_2}{p}|t|^{p}$ for any $t \in \mathbb{R}$. The rest of the proof is similar to (5), we omit it. \qed

\section{The Minimax and Compactness Lemmas}
In this section, we present a minimax lemma to obtain a bounded Palais-Smale sequence and a compactness lemma to verify the Palais-Smale condition.
We borrow some arguments of \cite{JeanjeanCV2020} and introduce the functional $\Phi:S(c)\rightarrow \mathbb{R} $ by
\begin{equation*}
  \Phi(u)= E_{\mu,c}\left(u_{t(u)}\right)
\end{equation*}
which satisfies
\begin{Lem}\label{TheAxuFunctional}
The functional $\Phi$ is of class $C^1$ and
$$ d \Phi(u)[\psi]=dE_{\mu,c}\left(u_{t(u)}\right)[\psi_{t(u)}]  $$
for any $u \in S(c)$ and $\psi \in T_{u}S(c)$, where $T_{u}S(c)$ the tangent space to $S(c)$ in $u$.
\end{Lem}
\Proof Similar to \cite[Lemma 5.4, Lemma 5.5, Lemma5.6]{Yao} and \cite[Lemma 4.2, Lemma 4.3]{JeanjeanCV2020}, we omit it.   \qed

We adopt a minimax argument as a variant of \cite[Theorem 4.5]{JeanjeanCV2020} or an application of \cite[Theorem 3.2]{Duality} to obtain a Palais-Smale sequence.
\begin{Def}[{\cite[Definition 3.1]{Duality}}]
Let $B$ be a closed subset of a metric space $X$. We say that a class $\mathcal{G}$ of compact subset of $X$ is a homotopy stable family with closed boundary $B$ provided
\begin{enumerate}
  \item [(i)] every set in $\mathcal{G}$ contains $B$.
  \item [(ii)] for any set $A \in \mathcal{G}$ and any homotopy $\eta \in C([0,1]\times X, X)$ that satisfies $\eta(t,u)=u$ for all $(t,u) \in \left(\{0\}\times X \right) \cup \left([0,1]\times B \right)$, one has $\eta(\{1\}\times A) \in \mathcal{G}$.
\end{enumerate}
\end{Def}
\begin{Prop}[{\cite[Theorem 3.2]{Duality}}]\label{Theorem32} Let $\varphi$ be a $C^1$-functional on a complete connected $C^1$-Finsler manifold $X$ (without boundary) and consider a homotopy-stable family $\mathcal{F}$ of compact subsets of $X$ with a closed boundary $B$.
Set $c=c(\varphi,\mathcal{F})=\inf\limits_{A \in \mathcal{F}}\max\limits_{x \in A}\varphi(x)$ and suppose that
$$ \sup\varphi(B)<c.  $$
Then, for every sequence of sets $(A_n)_n$ in $\mathcal{F}$ such that $\lim\limits_{n\rightarrow \infty}\sup\limits_{A_n}\varphi=c$, there exists a sequence $(x_n)_n$ in $X$ such that
\begin{enumerate}
  \item [(i)] $\lim\limits_{n \rightarrow \infty} \varphi(x_n)=c$
  \item [(ii)] $\lim\limits_{n \rightarrow \infty} \|d\varphi(x_n)\|=0$
  \item [(iii)] $\lim\limits_{n \rightarrow \infty} dist(x_n,A_n)=0$.
\end{enumerate}
Moreover, if $d\varphi$ is uniformly continuous, then $x_n$ can be chosen to be in $A_n$ for each $n$.
\end{Prop}
\noindent We remark that the case $B=\varnothing$ is admissible and $\sup\varphi(B)=-\infty$.
\begin{Lem}\label{ThePSsequence}
Let $N \geq 1$, $g$ satisfy (G1)--(G3), then
there exists a bounded Palais-Smale sequence $\{ u_n\} \subset S(c)$ for $E_{\mu,c}\big|_{S(c)}$ at the level $\E(\mu,c)$, and
  there exist a bounded sequence $\{\tilde{u}_n \}\subset \mathcal{P}(\mu,c)$ with $\tilde{u}_n(x)$ is a nonnegative radial and nonincreasing function of $|x|$ for every $n \geq 0$, such that $\| u_n -\tilde{u}_n \|_{H^1(\RN)} \rightarrow 0 $ as $n \rightarrow \infty$.
\end{Lem}
 \Proof Let $\mathcal{G}$ be the class of all singletons included in $S(c)$ and $B=\varnothing$. Then, there holds
 $$ \inf_{A\in\mathcal{G}}\max_{u\in A}\Phi(u) = \inf_{u \in S(c)} E_{\mu,c}\left(u_{t(u)}\right).      $$
For any $u \in S(c)$, we deduce from $u_{t(u)}\in \mathcal{P}(\mu,c)$ that $E_{\mu,c}\left(u_{t(u)}\right) \geq \E(\mu,c)$. Therefore, $\inf\limits_{A\in\mathcal{G}}\max\limits_{u\in A}\Phi(u) \geq \E(\mu,c)$.
On the other hand, for any $u \in \mathcal{P}(\mu,c)$, we have $t(u)=1$ and thus $E_{\mu,c}(u) \geq \inf\limits_{A\in\mathcal{G}}\max\limits_{u\in A}\Phi(u)$. This implies $\E(\mu,c) \geq \inf\limits_{A\in\mathcal{G}}\max\limits_{u\in A}\Phi(u)$.
Consequently, we obtain $$ \inf\limits_{A\in\mathcal{G}}\max\limits_{u\in A}\Phi(u) = \E(\mu,c)>-\infty. $$

Let $\{ v_n \}$ be an arbitrary minimizing sequence such that $\Phi(v_n)  \rightarrow \E(\mu,c)$, then $\{ (v_n)_{t(v_n)} \}$ is also a minimizing sequence of $\E(\mu,c)$. We may assume $\{ v_n\} \in \mathcal{P}(\mu,c)$ without loss of generality. Note that for any $n \geq 0$, let $|v_n|^*$ denote the Schwartz rearrangement of $|v_n|$, then
$$ \| |v_n|^* \|_2^2=\|v_n \|_2^2=c^2, \quad, \| \nabla |v_n|^* \|_2^2 \leq \| \nabla |v_n| \|_2^2 \leq \| \nabla v_n \|_2^2, $$
$$\int_{\RN}(I_{\alpha}*|v_n|^{\2})|v_n|^{\2}dx \leq \int_{\RN}(I_{\alpha}*||v_n|^*|^{\2})||v_n|^*|^{\2}dx,  $$
and for any $t > 0$,
$$ \int_{\RN}G((v_n)_t)dx=\int_{\RN}G\left((|v_n|^*)_t\right)dx.         $$
Therefore, we deduce for any $ t>0 $,
\begin{equation}\label{TheMax1}
   E_{\mu,c}\left((v_n)_{t}\right) \geq E_{\mu,c}\left(\left( |v_n|^*\right)_{t}\right).
\end{equation}
It follows from Lemma \ref{TheMaxFiber} that
there exists a unique $t(|v_n|^*)$ such that $\left( |v_n|^*\right)_{t(|v_n|^*)} \in \mathcal{P}(\mu,c)$. Let $\tilde{u}_n := \left( |v_n|^*\right)_{t(|v_n|^*)}$, then by Lemma \ref{TheMaxFiber} and (\ref{TheMax1}), one has
$$\aligned E_{\mu,c}(v_n)= \max_{t>0}E_{\mu,c}\left((v_n)_t\right) & \geq E_{\mu,c}\left((v_n)_{t(|v_n|^*)}\right) \\
& \geq E_{\mu,c}\left(\left( |v_n|^*\right)_{t(|v_n|^*)}\right) \geq \E(\mu,c). \endaligned         $$
Thus, by Lemma \ref{TheMaxFiber} (3), we obtain a bounded minimizing sequence $\{\tilde{u}_n \} \subset \mathcal{P}(\mu,c)$ of $\E(\mu,c)$ with $\tilde{u}_n(x)$ is radial and nonincreasing function of $|x|$ for every $n \geq 0$.

Now, using Proposition \ref{Theorem32}, we obtain a Palais-Smale sequence $\{u_n\} \subset S(c)$ for $\Phi$ at the level $\E(\mu,c)$ such that $dist_{H^1(\RN)}(u_n, \tilde{u}_n) \rightarrow 0$ as $n \rightarrow \infty$. It follows from Lemma \ref{TheAxuFunctional} that $\{ u_n\} \subset S(c)$ is
a bounded Palais-Smale sequence of $E_{\mu,c}\big|_{S(c)}$ at the level $\E(\mu,c)$.  \qed

Since $\{ u_n\}$ is bounded in $H^1(\RN)$, by Lemma \ref{ThePSsequence} and Lemma \ref{TheCompactnessTool}, there exists ${u}_{\mu,c} \in H^1(\RN)$, a radial and nonincreasing function of $|x|$, such that
$$\aligned &\tilde{u}_n \rightharpoonup {u}_{\mu,c} \text{\ in\ } H^1(\RN), \\
   & \tilde{u}_n \rightarrow {u}_{\mu,c}  \text{\ in\ } L^s(\RN)\text{\ for\ any\ }s\in(2,2^*), \\
   & \tilde{u}_n \rightarrow {u}_{\mu,c}\text{\ a.e.\ in\ } \RN.  \endaligned$$
Since $\| u_n -\tilde{u}_n \|_{H^1(\RN)} \rightarrow 0 $ as $n \rightarrow \infty$, we deduce that
$$\aligned & u_n \rightharpoonup {u}_{\mu,c}  \text{\ in\ } H^1(\RN), \\
& u_n \rightarrow {u}_{\mu,c}  \text{\ in\ } L^s(\RN)\text{\ for\ any\ }s\in(2,2^*),\\ & u_n \rightarrow {u}_{\mu,c} \text{\ a.e.\ in\ } \RN.  \endaligned$$
Thus, it follows that $P_{\mu,c}(u_n) \rightarrow 0$ as $n \rightarrow \infty$. Moreover, by the Lagrange multipliers rule, there exists $\{\lambda_n\} \subset \mathbb{R}$ such that for every $\omega \in H^1(\RN)$,
\begin{equation}\label{TheLagrangeMultipiers}
 \aligned &\int_{\RN}(\nabla u_n \nabla \omega + \lambda_n u_n \omega)dx -\mu \int_{\RN}g(u_n)\omega dx\\ &- \gamma \int_{\RN}(I_{\alpha}*|u_n|^{\2})|u_n|^{\2-2}u_n\omega dx = o(1)\|\omega \|.   \endaligned
\end{equation}
 Choose $\omega = u_n $ in (\ref{TheLagrangeMultipiers}), using \eqref{barGGg} and $P(u_n)\rightarrow0$, one has
\begin{equation}\label{TheLambda}\aligned
  \lambda_n c^2 &= \gamma \IIN + \mu \int_{\RN}g(u_n)u_n dx - \|\nabla u_n\|_2^2\\
                &\geq \mu \frac{q}{q-2} \int_{\RN}\bar{G}(u_n)dx - \|\nabla u_n\|_2^2 + o(1)\\
                &\geq \left(\frac{q}{q-2}\cdot\frac{2}{N}-1\right)\|\nabla u_n \|_2^2 + o(1).  \endaligned
\end{equation}
Since the boundedness of $\{u_n\}$ in $H^1(\RN)$ implies that $\lambda_n$ is bounded as well, we may assume $\lambda_n \rightarrow {\lambda}_{\mu,c} \geq 0$, up to a subsequence.

By the weak convergence, if $u_{\mu,c} \not\equiv 0$ and $\lambda_{\mu,c}>0$, then $u_{\mu,c}$ is a nontrivial solution to the fixed frequency problem
\begin{equation}\label{Thefixedfrequencygeneral}
   -\Delta {u} + {\lambda}_{\mu,c} {u} =\mu g(u) + \gamma (I_{\alpha}*|u|^{\2})| u|^{\2-2} u \quad \text{\ in\ }\RN.
\end{equation}
The corresponding action functional $A_{\mu,c}: H^1(\RN) \mapsto \mathbb{R}$ of (\ref{Thefixedfrequencygeneral}) is given by
\begin{equation}\label{theactionAgeneral}
\aligned  A_{\mu,c}(u) :=&\frac{1}{2}\|\nabla u\|_2^2 + \frac{1}{2}\lambda_{\mu,c}\| u\|_2^2 - \mu \int_{\RN}G(u)dx \\
                          & - \frac{N\gamma}{2(N+\alpha)} \II \\
                       =& E_{\mu,c}(u)+ \frac{1}{2}\lambda_{\mu,c}\| u\|_2^2. \endaligned
\end{equation}
It is well-known that all nontrivial solutions to (\ref{Thefixedfrequencygeneral}) belong to the Nehari manifold
\begin{equation}\label{Neharigeneral}
  \mathcal{N}(\mu,c):=\left\{ u \in H^1(\RN)\setminus\{0\}: N_{\mu,c}(u)=0   \right\}
\end{equation}
where
\begin{equation}\label{Nehariindentitygeneral}
   N_{\mu,c}(u):= \|\nabla u\|_2^2 + \lambda_{\mu,c}\|u \|_2^2 - \mu \int_{\RN}g(u)u dx - \gamma \II.
\end{equation}
\begin{Lem}\label{compactnesslemma}
The Palais-Smale sequence $\{u_n\}$ converges to $u_{\mu,c}$ strongly in $H^1(\RN)$ up to a subsequence, if
\begin{enumerate}
  \item [(1)] $u_{\mu,c} \not\equiv 0$ and $\lambda_{\mu,c}>0$.
  \item [(2)] $ \E(\mu,c)+\frac{1}{2}\lambda_{\mu,c}c^2 < a(\mu,c) + \frac{\alpha}{2(N+\alpha)} \gamma^{-\frac{N}{\alpha}}({\lambda}_{\mu,c} S_{\alpha})^{\frac{N+\alpha}{\alpha}}$ where
  $$a(\mu,c):= \inf_{u \in \mathcal{N}(\mu,c)} A_{\mu,c}(u). $$
\end{enumerate}
\end{Lem}
\Proof Note that $u_{\mu,c}$ is a nontrivial solution to (\ref{Thefixedfrequencygeneral}), one has $P_{\mu,c}(u_{\mu,c})=0$.  Let $w_n = u_n -  {u}_{\mu,c}$, then, by Brezis-Lieb lemma, there holds that
$$ \|\nabla u_n \|_2^2 =\|\nabla {u}_{\mu,c} \|_2^2 + \|\nabla w_n \|_2^2 + o(1), \quad  \| u_n \|_2^2 =\| {u}_{\mu,c} \|_2^2 + \| w_n \|_2^2 + o(1),    $$
$$\aligned  \IIN=&\int_{\RN}(I_{\alpha}*|w_n|^{\2})|w_n|^{\2}dx \\
                   &+ {\int_{\RN}(I_{\alpha}*|{u}_{\mu,c}|^{\2})|{u}_{\mu,c}|^{\2}dx}+o(1). \endaligned $$
Since $u_n \rightarrow {u}_{\mu,c}$ in $L^p(\RN)$ for $p \in (2,2^*)$ and $P_{\mu,c}(u_n)=P_{\mu,c}(u_{\mu,c})+P_{\mu,c}(w_n)+ o(1)$, we deduce $P_{\mu,c}(w_n)=o(1)$ and hence $\|\nabla w_n\|_2^2 =o(1)$. Using this fact and $N_{\mu,c}(u_{\mu,c})=0$, we obtain that
$$\aligned o(1) &= \|\nabla u_n\|_2^2 + \lambda_n \| u_n\|_2^2 - \gamma \IIN -\mu \int_{\RN}g(u_n)u_n dx  \\
                &= \lambda_n \|w_n \|_2^2 - \gamma \int_{\RN}(I_{\alpha}*|w_n|^{\2})|w_n|^{\2}dx. \endaligned$$
Now we may assume
$$ \theta = \lim_{n\rightarrow\infty} \lambda_n \| w_n \|_2^2 = \lim_{n\rightarrow\infty}\int_{\RN}(I_{\alpha}*|w_n|^{\2})|w_n|^{\2}dx, $$
and this implies $\theta=0$ or $\theta \geq \gamma^{-\frac{N}{\alpha}}({\lambda}_{\mu,c} S_{\alpha})^{\frac{N+\alpha}{\alpha}}$.
 If $\theta\geq \gamma^{-\frac{N}{\alpha}}({\lambda}_{\mu,c} S_{\alpha})^{\frac{N+\alpha}{\alpha}}$, then we have
\begin{equation*}
\aligned
  &\E(\mu,c)+\frac{{\lambda}_{\mu,c}}{2}c^2\\ &= E_{\mu,c}(u_n) + \frac{\lambda_n}{2}\|u_n\|_2^2 +o(1)\\
                                      &=A_{\mu,c}(u_{\mu,c})+\left( \frac{\lambda_n}{2}\|w_n\|_2^2 - \frac{N\gamma}{2(N+\alpha)}\int_{\RN}(I_{\alpha}*|w_n|^{\2})|w_n|^{\2}dx \right)+o(1)  \\
                                      &\geq a(\mu,c) +\frac{\alpha}{2(N+\alpha)} \gamma^{-\frac{N}{\alpha}}({\lambda}_{\mu,c} S_{\alpha})^{\frac{N+\alpha}{\alpha}} \endaligned
\end{equation*}
provided $A_{\mu,c}(u_{\mu,c}) \geq a(\mu,c)$ since $u_{\mu,c} \in \mathcal{N}_{\mu,c}$. However this contradicts (2) and we have $\theta=0$ which implies $u_n \rightarrow {u}_{\mu,c}$ strongly in $H^1(\RN)$ clearly.  \qed

\section{Existence of Ground states}
In this section, we give some estimates on $\E(\mu,c)$, $\lambda_{\mu,c}$ and $a(\mu,c)$ and obtain the normalized ground state solutions in two different cases.
\begin{Lem}\label{TheNableNormsmallPara}
Let $N \geq 1$, $g$ satisfy (G1)--(G4), then there exists a positive constant $\delta_1$ such that
\begin{equation}\label{fomsharpnablanormsmall}
\|\nabla u_n \|_2^2 \sim  \mu^{-\frac{4}{N(q-2)-4}} c^{-\frac{4N-2q(N-2)}{N(q-2)-4}}
\end{equation}
for $n$ sufficiently large and
\begin{equation}\label{fomsharpmmuc1}
  \E(\mu,c) + \frac{N\gamma}{2(N+\alpha)}c^{\frac{2(N+\alpha)}{N}}S_{\alpha}^{-\frac{N+\alpha}{N}} \sim \mu^{-\frac{4}{N(q-2)-4}} c^{-\frac{4N-2q(N-2)}{N(q-2)-4}}
\end{equation}
 if $\mathcal{N}_{p,q}(c,\mu)< \delta_1$.
\end{Lem}
\Proof It follows from (\ref{TheLowerConstant}) and (\ref{barGGg}) that
$$ \aligned      & \E(\mu,c) + \frac{N\gamma}{2(N+\alpha)}c^{\frac{2(N+\alpha)}{N}}S_{\alpha}^{-\frac{N+\alpha}{N}} +o(1) \\
            \geq & \frac{1}{2}\|\nabla u_n\|_2^2 - \mu \int_{\RN}G(u_n)dx\\
            \geq & \frac{N(p-2)-4}{2N(p-2)}\|\nabla u_n\|_2^2 \endaligned $$
and hence, we have $\|\nabla u_n\|_2^2 \leq C \mu^{-\frac{4}{N(q-2)-4}} c^{-\frac{4N-2q(N-2)}{N(q-2)-4}}  $ for $n$ sufficiently large by Lemma \ref{TheMaxFiber} (4).
On the other hand, from (G1)(G2) and  Lemma \ref{GNS}, we obtain that
\begin{equation}\label{muclargelowerbound}
\aligned  \|\nabla u\|_2^2 &= \mu \frac{N}{2}\int_{\RN}\bar{G}(u)dx \\
                           &\leq \mu C_{1} \left(c^{\frac{2q-N(q-2)}{2}}\|\nabla u\|_2^{\frac{N(q-2)}{2}}
                                             +c^{\frac{2p-N(p-2)}{2}}\|\nabla u\|_2^{\frac{N(p-2)}{2}} \right).   \endaligned
\end{equation}
Then there exists $\delta_1 >0$ such that
$$  C_1 C^{\frac{N(p-2)-4}{4}}  {\delta_1}^{\frac{N(p-2)-4}{4}} < \frac{1}{2},$$
and we deduce that
$\|\nabla u_n\|_2^2 \geq  C' \mu^{-\frac{4}{N(q-2)-4}} c^{-\frac{4N-2q(N-2)}{N(q-2)-4}}  $
 for $n$ sufficiently large by (\ref{muclargelowerbound}) if
$$ \mu^{1-\frac{N(p-2)-4}{N(q-2)-4}}c^{\frac{2N-p(N-2)}{2} -\frac{2N-q(N-2)}{N(q-2)-4}\frac{N(p-2)-4}{2}}
            = {\mathcal{N}_{p,q}(c,\mu)}^{\frac{N(p-2)-4}{4}}< {\delta_1}^{\frac{N(p-2)-4}{4}}. $$
Now we obtain (\ref{fomsharpnablanormsmall}).
Note that $E_{\mu,c}(u_n) \rightarrow \E(\mu,c)$ and $P_{\mu,c}(u_n) \rightarrow 0$ as $n \rightarrow +\infty$, it follows from Lemma \ref{GNS} and (\ref{TheLowerConstant}) that
$$\aligned &\E(\mu,c)+o_n(1)\\
                       =&\frac{1}{2}\|\nabla u_n \|_2^2 -\mu \int_{\RN}G(u_n)dx -\frac{\gamma N}{2(N+\alpha)}\IIN +o_n(1)\\
                       \geq&\left( \frac{1}{2}-\frac{2}{N(p-2)} \right)\|\nabla u_n\|_2^2 -\frac{\gamma N}{2(N+\alpha)}\IIN+o_n(1)\\
                       \geq& C \mu^{-\frac{4}{N(q-2)-4}} c^{-\frac{4N-2q(N-2)}{N(q-2)-4}}-\frac{N\gamma}{2(N+\alpha)}c^{\frac{2(N+\alpha)}{N}}S_{\alpha}^{-\frac{N+\alpha}{N}} \endaligned $$
provided (\ref{fomsharpnablanormsmall}). Combining this with Lemma \ref{TheMaxFiber} (4) we obtain (\ref{fomsharpmmuc1}).
This completes the proof. \qed

\begin{Lem}\label{convergesmall}
Let $N \geq 1$, $g$ satisfy (G1)--(G4), then there exists $\delta_2>0$ such that
$\{u_n \}$ converges to $u_{\mu,c}$ strongly in $H^1(\RN)$ up to a subsequence if
${\mathcal{N}_{p,q}(c,\mu)}< \delta_1$ and $\mathcal{M}_q(\mu,c)< \delta_2 $.
\end{Lem}
\Proof The proof is divided into three steps.

\textbf{Step 1.} We show ${\lambda}_{\mu,c}>0$ and $u_{\mu,c} \not\equiv 0$. If ${\lambda}_{\mu,c}=0$,
it follows from (\ref{TheLambda}) that $\|\nabla u_n\|_2^2=\|\nabla \tilde{u}_n \|_2^2=o_n(1)$ which contradicts with (\ref{fomsharpnablanormsmall}). Thus, we obtain ${\lambda}_{\mu,c}>0$.
If ${u}_{\mu,c}=0$, then from $\tilde{u}_n \rightarrow 0$ in $L^p(\RN)$ for any $p \in (2,2^*)$ and $P(\tilde{u}_n)=0$ we have $\|\nabla \tilde{u}_n \|_2^2 \rightarrow 0$ which also contradicts with (\ref{fomsharpnablanormsmall}). Thus, we obtain $u_{\mu,c} \not\equiv 0$.

\textbf{Step 2.} We estimate $ {\lambda}_{\mu,c}$ from above and below. It follows from (\ref{TheLambda}) and (\ref{fomsharpnablanormsmall}) that
\begin{equation*}
  {\lambda}_{\mu,c} \gtrsim  \mu^{-\frac{4}{N(q-2)-4}} c^{-\frac{4(q-2)}{N(q-2)-4}}.
\end{equation*}
Note that $E_{\mu,c}(u_n) \rightarrow \E(\mu,c)$ and $P_{\mu,c}(u_n) \rightarrow 0$, it follows from (G2) and (\ref{fomsharpnablanormsmall}) that
\begin{equation}\label{Theabove}
  \aligned \lambda_n c^2 +o(1) &=\gamma \IIN +\mu\int_{\RN}g(u_n)u_ndx-\|\nabla u_n\|_2^2\\
                         &\leq\gamma \IIN + \frac{2N-p(N-2)}{N(p-2)} \|\nabla u_n\|_2^2 \\
                         &\lesssim \gamma c^{\frac{2(N+\alpha)}{N}}S_{\alpha}^{-\frac{N+\alpha}{\alpha}} + C \mu^{-\frac{4}{N(q-2)-4}} c^{-\frac{4N-2q(N-2)}{N(q-2)-4}}  \\
                         &\lesssim \mu^{-\frac{4}{N(q-2)-4}} c^{-\frac{4N-2q(N-2)}{N(q-2)-4}}  \endaligned
\end{equation}
if $\mu^{\frac{2}{N(q-2)-4}} c^{\frac{\alpha}{N}+ \frac{2(q-2)}{N(q-2)-4}}< \delta_2' $ for some $\delta_2' >0$ suitably small.
Thus, we get
\begin{equation}\label{sharplambdasmall}
  \lambda_{\mu,c} c^2 \sim \mu^{-\frac{4}{N(q-2)-4}} c^{-\frac{4N-2q(N-2)}{N(q-2)-4}}.
\end{equation}

\textbf{Step 3.} We now show that $ u_n \rightarrow {u}_{\mu,c}$ strongly in $H^1(\RN)$, up to a subsequence. It follows from (\ref{fomsharpmmuc1}) and (\ref{sharplambdasmall}) that
\begin{equation*}
 \aligned
 \E(\mu,c)+ \frac{{\lambda}_{\mu,c}}{2}c^2 & \lesssim \mu^{-\frac{4}{N(q-2)-4}} c^{-\frac{4N-2q(N-2)}{N(q-2)-4}}-\frac{N\gamma}{2(N+\alpha)}c^{\frac{2(N+\alpha)}{N}}S_{\alpha}^{-\frac{N+\alpha}{\alpha}} \\
                                          & \lesssim \mu^{-\frac{4}{N(q-2)-4}} c^{-\frac{4N-2q(N-2)}{N(q-2)-4}}   \endaligned
\end{equation*}
and
$${\lambda}_{\mu,c}^{\frac{N+\alpha}{\alpha}} \gtrsim \left(\mu^{-\frac{4}{N(q-2)-4}} c^{-\frac{4(q-2)}{N(q-2)-4}} \right)^{\frac{N+\alpha}{\alpha}}. $$
Hence, there holds
\begin{equation}\label{TheKEY}
  \E(\mu,c)+ \frac{{\lambda}_{\mu,c}}{2}c^2 <  \frac{\alpha}{2(N+\alpha)} \gamma^{-\frac{N}{\alpha}}({\lambda}_{\mu,c} S_{\alpha})^{\frac{N+\alpha}{\alpha}}
\end{equation}
holds provided
$$   \mu^{\frac{4}{N(q-2)-4}\frac{N}{\alpha}}c^{\frac{4(q-2)}{N(q-2)-4}\frac{N+\alpha}{\alpha}-\frac{4q-2N(q-2)}{N(q-2)-4}}
            =  \mathcal{M}_{q}(c,\mu))^{\frac{N}{\alpha}} <{\delta_2}^{\frac{N}{\alpha}}  $$
for some suitable $\delta_2>0$. Since it is standard to see $a(\mu,c) \geq 0$ for any $\mu>0$ and $c>0$, it follows from (\ref{TheKEY}) that
 $$ \E(\mu,c)+\frac{1}{2}\lambda_{\mu,c}c^2 < a(\mu,c) + \frac{\alpha}{2(N+\alpha)} \gamma^{-\frac{N}{\alpha}}({\lambda}_{\mu,c} S_{\alpha})^{\frac{N+\alpha}{\alpha}}. $$
By Lemma \ref{compactnesslemma}, we deduce that $u_n \rightarrow {u}_{\mu,c}$ strongly in $H^1(\RN)$. This completes the proof.
\qed

\noindent\textbf{The proof of Theorem \ref{TheSecondResult}.}
Let $ \delta=\min\{\delta_1, \delta_2\}$ and by Lemma \ref{convergesmall}, we obtain a normalized solution $({u}_{\mu,c}, {\lambda}_{\mu,c}) \in H^1(\RN) \times \mathbb{R}^+$ satisfying
$$ E_{\mu,c}({u}_{\mu,c})=\E(\mu,c)= \inf\{  E_{\mu,c}(v): v\in S(c), (E_{\mu,c}|_{S(c)})'(v)=0    \}. $$
Recall that $\tilde{u}_n \rightarrow {u}_{\mu,c}$ a.e. in $\RN$ and $\tilde{u}_n \geq 0$ for any $n \geq 0$, we deduce ${u}_{\mu,c} \geq 0$ and by the maximum principle $ {u}_{\mu,c} > 0 $. Then
${u}_{\mu,c}$ is a positive, radial, nonincreasing function of $|x|$.
The estimate (\ref{TheCharOfLag}) comes from (\ref{sharplambdasmall}), and (\ref{TheCharOfEnergy}) comes from (\ref{fomsharpmmuc1}). \qed

\begin{Lem}\label{TheNableNormlargePara}
Let $N \geq 1$, $g$ satisfy (G1)--(G3) and (G5), then there exists a positive constant $\tau_1$ such that
\begin{equation}\label{fomsharpnablanormlarge}
\|\nabla u_n \|_2^2 \sim  \mu^{-\frac{4}{N(p-2)-4}} c^{-\frac{4N-2p(N-2)}{N(p-2)-4}}
\end{equation}
for $n$ sufficiently large and
\begin{equation}\label{fomsharpmmuc1large}
  \E(\mu,c) + \frac{N\gamma}{2(N+\alpha)}c^{\frac{2(N+\alpha)}{N}}S_{\alpha}^{-\frac{N+\alpha}{N}} \sim \mu^{-\frac{4}{N(p-2)-4}} c^{-\frac{4N-2p(N-2)}{N(p-2)-4}}
\end{equation}
 if $ {\mathcal{N}_{p,q}(c,\mu)} > \tau_1$.
\end{Lem}
\Proof Argue as Lemma \ref{TheNableNormsmallPara} with Lemma \ref{TheMaxFiber} (5) in hand, we omit it.   \qed

\begin{Lem}\label{convergelarge}
Let $N \geq 1$, $g$ satisfy (G1)--(G3) and (G5),  then
$\{u_n \}$ converges to $u_{\mu,c}$ strongly in $H^1(\RN)$ up to a subsequence if  $ {\mathcal{N}_{p,q}(c,\mu)} > \tau_1$ and
$$ \min\left\{  \mathcal{M}_q(c,\mu), \mathcal{M}_p(c,\mu)\right\}> \tau_2. $$
\end{Lem}
\Proof Similar to Lemma \ref{convergesmall} we have $u_{\mu,c}\not\equiv0$ and $\lambda_{\mu,c}>0$, so it is sufficient to verify the inequality
$$ \E(\mu,c)+\frac{1}{2}\lambda_{\mu,c}c^2 < a(\mu,c) + \frac{\alpha}{2(N+\alpha)} \gamma^{-\frac{N}{\alpha}}({\lambda}_{\mu,c} S_{\alpha})^{\frac{N+\alpha}{\alpha}}. $$

\textbf{Step 1.} We estimate $ {\lambda}_{\mu,c}$ from above and below.
Since $$ \lambda_{\mu,c}c^2 + o(1) = \mu\int_{\RN}g(u_n)u_ndx +\gamma \IIN -\|\nabla u_n \|_2^2, $$
it follows from (G2), (\ref{fomsharpnablanormlarge}) and (\ref{TheLowerConstant}) that
$$  \lambda_{\mu,c}c^2 -\gamma S_{\alpha}^{-\frac{N+\alpha}{N}}c^{\frac{2(N+\alpha)}{N}} \leq C_1 \mu^{-\frac{4}{N(p-2)-4}} c^{-\frac{4N-2p(N-2)}{N(p-2)-4}} $$
for some $C_1>0$. On the other hand, notice that
$$ \E(\mu,c)+o(1) = \frac{1}{2}\|\nabla u\|_2^2 - \mu \int_{\RN}G(u_n)dx - \frac{N\gamma}{2(N+\alpha)}\IIN, $$
we deduce that
\begin{equation}\label{lambdamularge}
 \aligned & \lambda_{\mu,c}c^2 + \frac{2(N+\alpha)}{N}\E(\mu,c)+o(1) \\
         =&  \mu \int_{\RN}g(u_n)u_n dx - \mu \frac{2(N+\alpha)}{N}\int_{\RN}G(u_n)dx +\frac{N+\alpha}{N}\|\nabla u_n\|_2^2 -\|\nabla u_n\|_2^2 +o(1) \\
         \geq & \left[ \frac{2N-q(N-2)}{N(q-2)}+\frac{N+\alpha}{N}\frac{N(p-2)-4}{2N(p-2)}  \right]\|\nabla u_n\|_2^2 +o(1)\\
         \gtrsim & \mu^{-\frac{4}{N(p-2)-4}} c^{-\frac{4N-2p(N-2)}{N(p-2)-4}}. \endaligned
\end{equation}
 if $ (\mu^N c^{4})^{q-p}> \tau_1$.
Combining (\ref{lambdamularge}) with (\ref{fomsharpmmuc1large}) one has
 $$   \lambda_{\mu,c}c^2 -\gamma S_{\alpha}^{-\frac{N+\alpha}{N}}c^{\frac{2(N+\alpha)}{N}} \geq -C_2 \mu^{-\frac{4}{N(p-2)-4}} c^{-\frac{4N-2p(N-2)}{N(p-2)-4}}  $$
for some $C_2 > 0$. Thus, we deduce that
\begin{equation}\label{largeLagrange}
   \Big|\lambda_{\mu,c}c^2 -\gamma S_{\alpha}^{-\frac{N+\alpha}{N}}c^{\frac{2(N+\alpha)}{N}}\Big| \leq O \left(\mu^{-\frac{4}{N(p-2)-4}} c^{-\frac{4N-2p(N-2)}{N(p-2)-4}}\right).
\end{equation}
Then there exists a constant $\tau_2' >0$ such that
\begin{equation}\label{lambdac}
   \frac{1}{2} \gamma S_{\alpha}^{-\frac{N+\alpha}{N}}c^{\frac{2\alpha}{N}} \leq \lambda_{\mu,c} \leq \frac{3}{2} \gamma S_{\alpha}^{-\frac{N+\alpha}{N}} c^{\frac{2\alpha}{N}}
\end{equation}
if
\begin{equation}\label{lowerbound1}
   \mu^{\frac{4}{N(p-2)-4}} c^{\frac{2\alpha}{N}+ \frac{4(p-2)}{N(p-2)-4}} > \tau_2'.
\end{equation}

\textbf{Step 2.}  We estimate $ a({\mu,c})$ from below. By the definition of $a(\mu,c)$ let $\{ v_{\mu,c} \} \subset \mathcal{N}_{\mu,c}$ be a family of functions such that
$$  a(\mu,c) \leq A_{\mu,c}(v_{\mu,c}) \leq \frac{3}{2}a(\mu,c).  $$
Then it is standard to see $ a(\mu,c) \gtrsim \|\nabla v_{\mu,c} \|_2^2 +\lambda_{\mu,c}\|v_{\mu,c}\|_2^2$.

Case i)\ $\lambda_{\mu,c} \mu^{-1} \leq 1$. In this case, let $\tilde{v}_{\mu,c} = \mu^{\frac{1}{p-2}}\lambda_{\mu,c}^{-\frac{1}{p-2}} v( \lambda_{\mu,c}^{-\frac{1}{2}}x)$, then $\tilde{v}_{\mu,c}$ satisfies
\begin{equation}\label{lowerboundscaling}
\|\nabla v_{\mu,c} \|_2^2 +\lambda_{\mu,c}\|v_{\mu,c}\|_2^2=\mu^{-\frac{2}{p-2}}\lambda_{\mu,c}^{\frac{2N-p(N-2)}{2(p-2)}}  \Big(\|\nabla \tilde{v}_{\mu,c}\|_2^2 + \|\tilde{v}_{\mu,c} \|_2^2\Big).
\end{equation}
It follows from (\ref{lambdac}) that there exists a constant $\tau_3'>0$ such that
$$  \mu^{-\frac{2\alpha}{N(p-2)}} \lambda_{\mu,c}^{\frac{2\alpha}{N(p-2)}-1-\frac{\alpha}{2}} \leq 1 $$
if
\begin{equation}\label{lowerbound2}
   \mu^{\frac{4}{N(p-2)-4}} c^{\frac{2\alpha}{N}+ \frac{4(p-2)}{N(p-2)-4}} > \tau_3'.
\end{equation}
Thus, one has
$$ \aligned \|\nabla \tilde{v}_{\mu,c}\|_2^2 + \|\tilde{v}_{\mu,c} \|_2^2  = &\lambda_{\mu,c}^{-\frac{p-1}{p-2}}\mu^{\frac{p-1}{p-2}}\int_{\RN}g(\lambda_{\mu,c}^{\frac{1}{p-2}}\mu^{-\frac{1}{p-2}}\tilde{v}_{\mu,c})\lambda_{\mu,c}^{\frac{1}{p-2}}\mu^{-\frac{1}{p-2}}\tilde{v}_{\mu,c} dx \\ &  + \gamma\mu^{-\frac{2\alpha}{N(p-2)}} \lambda_{\mu,c}^{\frac{2\alpha}{N(p-2)}-1-\frac{\alpha}{2}}\int_{\RN}(I_{\alpha}*|\tilde{v}_{\mu,c}|^{\2})|\tilde{v}_{\mu,c}|^{\2}dx \\
\leq & \|\tilde{v}_{\mu,c}\|_p^p + \|\tilde{v}_{\mu,c}\|_q^q + \int_{\RN}(I_{\alpha}*|\tilde{v}_{\mu,c}|^{\2})|\tilde{v}_{\mu,c}|^{\2}dx\\
\lesssim &\left(\|\tilde{v}_{\mu,c} \|_{H^1}^p + \|\tilde{v}_{\mu,c} \|_{H^1}^q+\|\tilde{v}_{\mu,c} \|_{H^1}^{\frac{2(N+\alpha)}{N}}   \right), \endaligned $$
and hence $ \|\nabla \tilde{v}_{\mu,c}\|_2^2 + \|v_{\mu,c} \|_2^2 \gtrsim 1$ since $q>p>2$. From (\ref{lambdac}) and (\ref{lowerboundscaling}) we obtain
\begin{equation}\label{loweramc1}
  a(\mu,c) \gtrsim \mu^{-\frac{2}{p-2}}c^{\frac{\alpha}{N}\frac{2N-p(N-2)}{p-2}}.
\end{equation}

Case ii)\ $\lambda_{\mu,c}^{-1} \mu <1$.  In this case, let $\bar{v}_{\mu,c} = \mu^{\frac{1}{q-2}}\lambda_{\mu,c}^{-\frac{1}{q-2}} v( \lambda_{\mu,c}^{-\frac{1}{2}}x)$, then $\bar{v}_{\mu,c}$ satisfies
\begin{equation}\label{lowerboundscaling2}
\|\nabla v_{\mu,c} \|_2^2 +\lambda_{\mu,c}\|v_{\mu,c}\|_2^2=\mu^{-\frac{2}{q-2}}\lambda_{\mu,c}^{\frac{2N-q(N-2)}{2(q-2)}}  \Big(\|\nabla \bar{v}_{\mu,c}\|_2^2 + \|\bar{v}_{\mu,c} \|_2^2\Big).
\end{equation}
It follows from (\ref{lambdac}) that there exists a constant $\tau_3''>0$ such that
$$  \mu^{-\frac{2\alpha}{N(q-2)}} \lambda_{\mu,c}^{\frac{2\alpha}{N(q-2)}-1-\frac{\alpha}{2}} \leq 1 $$
if
\begin{equation}\label{lowerbound22}
   \mu^{\frac{4}{N(q-2)-4}} c^{\frac{2\alpha}{N}+ \frac{4(q-2)}{N(q-2)-4}} > \tau_3''.
\end{equation}
Thus, one has
$$ \aligned \|\nabla \bar{v}_{\mu,c}\|_2^2 + \|\bar{v}_{\mu,c} \|_2^2  = &\lambda_{\mu,c}^{-\frac{q-1}{q-2}}\mu^{\frac{q-1}{q-2}}g(\lambda_{\mu,c}^{\frac{1}{q-2}}\mu^{-\frac{1}{q-2}}\bar{v}_{\mu,c}) \\ &  + \gamma\mu^{-\frac{2\alpha}{N(q-2)}} \lambda_{\mu,c}^{\frac{2\alpha}{N(q-2)}-1-\frac{\alpha}{2}}\int_{\RN}(I_{\alpha}*|\bar{v}_{\mu,c}|^{\2})|\bar{v}_{\mu,c}|^{\2}dx \\
\leq & \|\bar{v}_{\mu,c}\|_p^p + \|\bar{v}_{\mu,c}\|_q^q + \int_{\RN}(I_{\alpha}*|\bar{v}_{\mu,c}|^{\2})|\bar{v}_{\mu,c}|^{\2}dx\\
\lesssim  &\left(\|\tilde{v}_{\mu,c} \|_{H^1}^p + \|\tilde{v}_{\mu,c} \|_{H^1}^q+\|\tilde{v}_{\mu,c} \|_{H^1}^{\frac{2(N+\alpha)}{N}}   \right), \endaligned $$
and hence $ \|\nabla \bar{v}_{\mu,c}\|_2^2 + \|\bar{v}_{\mu,c} \|_2^2 \gtrsim 1$. From (\ref{lambdac}) and (\ref{lowerboundscaling2}) we obtain
\begin{equation}\label{loweramc2}
  a(\mu,c) \gtrsim \mu^{-\frac{2}{q-2}}c^{\frac{\alpha}{N}\frac{2N-q(N-2)}{q-2}}.
\end{equation}

\textbf{Step 3.} From (\ref{fomsharpmmuc1large}) and (\ref{largeLagrange}) we obtain
\begin{equation}\label{thestrict1}
 \aligned  \E(\mu,c)+\frac{1}{2}\lambda_{\mu,c}c^2 \leq C \mu^{-\frac{4}{N(p-2)-4}} c^{-\frac{4N-2p(N-2)}{N(p-2)-4}} +\gamma \frac{\alpha}{2(N+\alpha)} S_{\alpha}^{-\frac{N+\alpha}{N}}c^{\frac{2(N+\alpha)}{N}} \endaligned
\end{equation}
and
\begin{equation}\label{thestrict2}
\aligned  & a(\mu,c) + \frac{\alpha}{2(N+\alpha)} \gamma^{-\frac{N}{\alpha}}({\lambda}_{\mu,c} S_{\alpha})^{\frac{N+\alpha}{\alpha}} \\
 \geq  &a(\mu,c)-C' \mu^{-\frac{4}{N(p-2)-4}} c^{-\frac{4N-2p(N-2)}{N(p-2)-4}} + \gamma\frac{\alpha}{2(N+\alpha)} S_{\alpha}^{-\frac{N+\alpha}{N}}c^{\frac{2(N+\alpha)}{N}}.   \endaligned
\end{equation}
As a consequence of (\ref{loweramc1}) and (\ref{loweramc2}) we deduce that
\begin{equation}\label{process}
   \E(\mu,c)+\frac{1}{2}\lambda_{\mu,c}c^2 < a(\mu,c) + \frac{\alpha}{2(N+\alpha)} \gamma^{-\frac{N}{\alpha}}({\lambda}_{\mu,c} S_{\alpha})^{\frac{N+\alpha}{\alpha}}
\end{equation}
if $a(\mu,c) \gtrsim \min\{ \mu^{-\frac{2}{q-2}}c^{\frac{\alpha}{N}\frac{2N-q(N-2)}{q-2}},\mu^{-\frac{2}{p-2}}c^{\frac{\alpha}{N}\frac{2N-p(N-2)}{p-2}} \} \gtrsim  \mu^{-\frac{4}{N(p-2)-4}} c^{-\frac{4N-2p(N-2)}{N(p-2)-4}}$. Since
$$\aligned & \mu^{\frac{4}{N(p-2)-4}} c^{\frac{4N-2p(N-2)}{N(p-2)-4}} \mu^{-\frac{2}{p-2}}c^{\frac{\alpha}{N}\frac{2N-p(N-2)}{p-2}} \\
           =&  \left( \mathcal{M}_p(c,\mu)   \right)^{\frac{2N-p(N-2)}{2(p-2)}}> \left( \min\{ \tau_2',\tau_3' \}\right)^{\frac{2N-p(N-2)}{2(p-2)}} \endaligned $$
and
$$ \aligned &\mu^{\frac{4}{N(p-2)-4}} c^{\frac{4N-2p(N-2)}{N(p-2)-4}} \mu^{-\frac{2}{q-2}}c^{\frac{\alpha}{N}\frac{2N-q(N-2)}{q-2}}
=\mathcal{N}_{p,q}(c,\mu) \left( \mathcal{M}_q(c,\mu)    \right)^{\frac{2N-q(N-2)}{2(q-2)}}
>\tau_1  \left( \tau_3'' \right)^{\frac{2N-q(N-2)}{2(q-2)}},  \endaligned $$
take $\tau_2', \tau_3', \tau_3'' $ suitably large and let $\tau_2=\max\{ \tau_2', \tau_3', \tau_3'' \}$, we deduce (\ref{process}).
 By Lemma \ref{compactnesslemma}, we obtain that $u_n \rightarrow {u}_{\mu,c}$ strongly in $H^1(\RN)$ and this completes the proof. \qed

\noindent
\textbf{The proof of Theorem \ref{theresult2}.} Let $\tau=\max\{ 1, \tau_1, \tau_2 \}$, then $\mathcal{M}_p(c,\mu)=\mathcal{N}_{p,q}(c,\mu)\cdot\mathcal{M}_q(c,\mu)>\tau^2>\tau_2$. By Lemma \ref{convergelarge}, we obtain a normalized ground state solution $({u}_{\mu,c}, {\lambda}_{\mu,c}) \in H^1(\RN) \times \mathbb{R}^+$, and ${u}_{\mu,c}$ is a positive, radial, nonincreasing function of $|x|$. The estimates (\ref{TheCharLag2}) and (\ref{TheCharOfEnergy2}) come from (\ref{fomsharpnablanormlarge}) and (\ref{fomsharpmmuc1large}).   \qed

\section{The asymptotic profiles}
In this section, we establish the asymptotic behavior of the solutions $\{u_{\mu,c}\}$ as $c$ and $\mu$ varies.
Firstly we introduce the rescaling
\begin{equation}\label{Thevscaling1}
  v:=\mu^{\frac{N}{N(q-2)-4}}c^{\frac{4}{N(q-2)-4}}u(\mu^{\frac{2}{N(q-2)-4}}c^{\frac{2(q-2)}{N(q-2)-4}} x).
\end{equation}
Then the mass-constraint problem (\ref{ThePro}) is reduced to
\begin{equation}\label{LIMITPro}
    \left\{ \aligned  &-\Delta v + \lambda v = g_{\mu,c}(v) +\gamma \mathcal{M}_q(c,\mu) (I_\alpha * |v|^{\frac{N+\alpha}{N}})|v|^{\frac{N+\alpha}{N}-2}v   & \text{in\ \ }  \mathbb{R}^N,\\
    & \int_{\mathbb{R}^N}v^2 dx =1, \quad v \in H^1(\mathbb{R}^N), \endaligned \right.
\end{equation}
where $g_{\mu,c}(v)=\mu^{\frac{N(q-1)}{N(q-2)-4}}c^{\frac{4(q-1)}{N(q-2)-4}}g(\mu^{-\frac{N}{N(q-2)-4}}c^{-\frac{4}{N(q-2)-4}}v)$. The corresponding energy functional is defined by
\begin{equation}\label{LIMITfunctional}
\aligned  J_{\mu,c}(v)=& \frac{1}{2}\|\nabla v\|_2^2 - \int_{\RN}G_{\mu,c}(v)dx \\
                       &- \frac{N\gamma}{2(N+\alpha)} \mathcal{M}_q(c,\mu) \int_{\RN}(I_\alpha * |v|^{\frac{N+\alpha}{N}})|v|^{\frac{N+\alpha}{N}}dx  \endaligned
\end{equation}
where $$ G_{\mu,c}(v)=\mu^{\frac{Nq}{N(q-2)-4}}c^{\frac{4q}{N(q-2)-4}}G(\mu^{-\frac{N}{N(q-2)-4}}c^{-\frac{4}{N(q-2)-4}}v) $$
respectively. Then $J_{\mu,c}(v)=\mu^{\frac{4}{N(q-2)-4}} c^{\frac{4q-2N(q-2)}{N(q-2)-4}}E_{\mu,c}(u)$. The Poho\v{z}aev manifold is defined by
\begin{equation}\label{LIMITmanifold}
  \tilde{\mathcal{P}}(\mu,c)=\{u \in S(1):\tilde{P}_{\mu,c}(v)=0 \}
\end{equation}
where
\begin{equation}\label{LIMITPoho}
 \aligned \tilde{P}_{\mu,c}(v)&=\|\nabla v\|_2^2 - \mu^{\frac{Nq}{N(q-2)-4}}c^{\frac{4q}{N(q-2)-4}}\frac{N}{2}\int_{\RN}\bar{G}(\mu^{-\frac{N}{N(q-2)-4}}c^{-\frac{4}{N(q-2)-4}}v)dx\\
                              &=\|\nabla v\|_2^2 - \frac{N}{2}\int_{\RN}\bar{G}_{\mu,c}(v)dx.    \endaligned
\end{equation}

\begin{Lem}\label{calcu5.1}
Let $\{({u}_{\mu,c}, {\lambda}_{\mu,c})\}$ denotes a family of the normalized ground state solutions to (\ref{ThePro}) under the assumptions of Theorem \ref{TheSecondResult}, and $v_{\mu,c}$ is the rescaling (\ref{Thevscaling1}) of ${u}_{\mu,c}$. Then for $\mathcal{M}_q(c,\mu)$ sufficiently small
\begin{enumerate}
  \item [(1)] $\|v_{\mu,c} \|_2^2=c^{-2}\|{u}_{\mu,c}\|_2^2 =1$, $\| v_{\mu,c}  \|_q^q =\mu^{\frac{4}{N(q-2)-4}} c^{\frac{4q-2N(q-2)}{N(q-2)-4}}\| {u}_{\mu,c} \|_q^q$;
  \item [(2)] $\| \nabla v_{\mu,c}  \|_2^2 =\mu^{\frac{4}{N(q-2)-4}} c^{\frac{4q-2N(q-2)}{N(q-2)-4}} \|\nabla {u}_{\mu,c} \|_2^2$;
  \item [(3)] $\int_{\RN}(I_{\alpha}*|v_{\mu,c}|^{\2})|v_{\mu,c}|^{\2}dx = c^{-\frac{2(N+\alpha)}{N}} \int_{\RN}(I_{\alpha}*|{u}_{\mu,c}|^{\2})|{u}_{\mu,c}|^{\2}dx$;
  \item [(4)] $J_{\mu,c}(v_{\mu,c})= \inf_{v \in   \tilde{\mathcal{P}}(\mu,c)} J_{\mu,c}(v)= \mu^{\frac{4}{N(q-2)-4}} c^{\frac{4q-2N(q-2)}{N(q-2)-4}} \E(\mu,c)$.
\end{enumerate}
\end{Lem}
\Proof (1)(2)(3) come from direct calculation. Note that $v \in \tilde{\mathcal{P}}(\mu,c)$ if and only if $u \in \mathcal{P}(\mu,c)$ and
$J_{\mu,c}(v_{\mu,c})=\mu^{\frac{4}{N(q-2)-4}} c^{\frac{4q-2N(q-2)}{N(q-2)-4}}E_{\mu,c}({u}_{\mu,c})$, then (4) follows. This completes the proof. \qed
\begin{Lem}\label{aspG}
Let $\mathcal{N}_{p,q}(c,\mu) \rightarrow 0$, then
\begin{equation}\label{LIMITG}
  \int_{\RN}G_{\mu,c}(v_{\mu,c})dx=\frac{1}{q}\|v_{\mu,c}\|_q^q+o_{\mu,c}(1)
\end{equation}
and
\begin{equation}\label{LIMITBarG}
 \frac{N}{2}\int_{\RN}\bar{G}_{\mu,c}(v_{\mu,c})dx=\frac{N(q-2)}{2q}\|v_{\mu,c}\|_q^q+o_{\mu,c}(1)
\end{equation}
where $o_{\mu,c}(1)\rightarrow 0$ as $\mathcal{N}_{p,q}(c,\mu) \rightarrow 0$.
\end{Lem}
\Proof Let $g(t):=t^{q-1}+\tilde{f}(t)$ for $t \geq 0$, then $G(t)=\frac{1}{q}t^q + \tilde{F}(t)$, and
$$ \lim_{t \rightarrow +\infty}\tilde{f}(t)/|t|^{q-1}=  \lim_{t \rightarrow +\infty}\tilde{F}(t)/|t|^q= 0. $$
Therefore, for any $\varepsilon>0$, there exists a constant $C_{\varepsilon}>0$ such that
$$ \Big|G(t)-\frac{1}{q}|t|^q \Big| \leq |\tilde{F}(t)| \leq \varepsilon|t|^{q} + C_{\varepsilon}|t|^{2+\frac{4}{N}}, $$
and hence
\begin{equation}\label{EstiG}
\aligned \Bigg| \int_{\RN}G_{\mu,c}(v_{\mu,c})dx-\frac{1}{q}\|v_{\mu,c}\|_q^q \Bigg|\leq & \varepsilon\|v_{\mu,c}\|_q^q + \mu c^{\frac{4}{N}}C_{\varepsilon}\|v_{\mu,c}\|_{2+\frac{4}{N}}^{2+\frac{4}{N}}. \endaligned
\end{equation}
By (\ref{TheCharOfEnergy}) and Lemma \ref{calcu5.1} (1)(2), one has $\{v_{\mu,c}\}$ is bounded in $H^1(\RN)$. Clearly, due to the arbitrariness of $\varepsilon$ and the boundedness of $\{v_{\mu,c}\}$,
(\ref{EstiG}) implies
$$  \int_{\RN}G_{\mu,c}(v_{\mu,c})dx=\frac{1}{q}\|v_{\mu,c}\|_q^q+o_{\mu,c}(1) .$$
Similarly, for any $\varepsilon>0$, there exists a constant $C_{\varepsilon}>0$ such that
$$ \Big|\frac{N}{2}\bar{G}(t)-\frac{N(q-2)}{2q}|t|^q \Big| \leq \Big|\frac{N}{2}\tilde{f}(t)t-N\tilde{F}(t)\Big| \leq \varepsilon|t|^{q} + C_{\varepsilon}|t|^{2+\frac{4}{N}}, $$
and hence
$$\frac{N}{2}\int_{\RN}\bar{G}_{\mu,c}(v_{\mu,c})dx=\frac{N(q-2)}{2q}\|v_{\mu,c}\|_q^q+o_{\mu,c}(1). $$
This completes the proof. \qed

Recall that $e(1,1)$ admits a unique minimizer, that is, the unique solution $(Z, \Lambda)$ to
\begin{equation}\label{ThelimitPro}
    \left\{ \aligned  &-\Delta u + \lambda u = |u|^{q-2}u \quad \text{in\ \ }  \mathbb{R}^N,\\
    & \int_{\mathbb{R}^N}u^2 dx =1, \quad u \in H^1(\mathbb{R}^N), \endaligned \right.
    \end{equation}
and $I(Z)=\inf_{u \in \mathcal{M}}I(u)=e(1,1)$
where
\begin{equation}\label{Zfunctional}
  I(u)=\frac{1}{2}\|\nabla u\|_2^2 -\frac{1}{q}\|u\|_q^q
\end{equation}
and $$\mathcal{M}=\left\{ u \in S(c): \|\nabla u\|_2^2=\frac{N(q-2)}{2q}\|u\|_q^q \right\}.$$
\begin{Lem}\label{Lemma5.3}
Let $\mathcal{N}_{p,q}(c,\mu) \rightarrow 0$, then there exists a unique $t_{\mu,c}>0$ such that $(v_{\mu,c})_{t_{\mu,c}} \in \mathcal{M}$ and $t_{\mu,c} \rightarrow 1$ as
 $\mathcal{N}_{p,q}(c,\mu) \rightarrow 0$.
\end{Lem}
\Proof
It is easy to see $$t_{\mu,c}=\left( \frac{2q}{N(q-2)}\frac{\|\nabla v_{\mu,c}\|_2^2}{\|v_{\mu,c}\|_q^q} \right)^{\frac{2}{N(q-2)-4}}.    $$
Since $v_{\mu,c} \in \tilde{\mathcal{P}}(\mu,c)$, that is,
 $$ \|\nabla v_{\mu,c}\|_2^2 = \frac{N}{2}\int_{\RN}\bar{G}_{\mu,c}(v_{\mu,c})dx,$$
and $\|v_{\mu,c}\|_q^q \gtrsim 1$, it follows from (\ref{LIMITBarG}) that
$$t_{\mu,c}=\left( \frac{2q}{N(q-2)}\frac{\frac{N}{2}\int_{\RN}\bar{G}_{\mu,c}(v_{\mu,c})dx}{\|v_{\mu,c}\|_q^q} \right)^{\frac{2}{N(q-2)-4}}\rightarrow 1    $$
as $\mathcal{N}_{p,q}(c,\mu) \rightarrow 0$.    \qed
\begin{Lem}\label{Lemma5.4}
Let $\mathcal{N}_{p,q}(c,\mu) \rightarrow 0$, then there exists a unique $s_{\mu,c}>0$ such that $(Z)_{s_{\mu,c}} \in \tilde{\mathcal{P}}(\mu,c)$ and $s_{\mu,c} \rightarrow 1$ as
 $\mathcal{N}_{p,q}(c,\mu) \rightarrow 0$.
\end{Lem}
\Proof The existence and uniqueness of $s_{\mu,c}$ come from Lemma \ref{TheMaxFiber}. Moreover, by (G1)(G2)(G4), there exists a constant $C>0$ such that
$$ \frac{N}{2}\bar{G}(t) \geq \frac{N(p-2)}{2}G(s)\geq C|t|^q \quad \text{\ for\ all\ } t \in \mathbb{R}.$$
Then, we deduce
\begin{equation}\label{sbound}
  s_{\mu,c}^2\|\nabla Z\|_2^2=\frac{N}{2}\int_{\RN}\bar{G}_{\mu,c}\left((Z)_{s_{\mu,c}}\right)dx \geq Cs_{\mu,c}^{\frac{N(q-2)}{2}}\|Z\|_q^q
\end{equation}
and $\{s_{\mu,c}\}$ is bounded in $\mathbb{R}$. By (G1)(G2) we deduce that $s_{\mu,c} \gtrsim 1$ for $\mathcal{N}_{p,q}(c,\mu)$ sufficiently small.
Now we may assume $s_{\mu,c} \rightarrow s_0$ as $\mathcal{N}_{p,q}(c,\mu)\rightarrow 0$, up to a subsequence, then
we obtain that
$$s_0^2 \|\nabla Z\|_2^2 = s_0^{\frac{N(q-2)}{2}}\frac{N(q-2)}{2q}\|Z\|_q^q  $$
by letting $\mathcal{N}_{p,q}(c,\mu) \rightarrow 0$ and adopting the similar argument in Lemma \ref{aspG}. Note that $Z$ is the unique solution to (\ref{ThelimitPro}), one has
$\|\nabla Z\|_2^2=\frac{N(q-2)}{2q}\| Z\|_q^q$ and it follows that $s_0=1$. Therefore, we obtain
$$\liminf\limits_{\mathcal{N}_{p,q}(c,\mu) \rightarrow 0}s_{\mu,c}=\limsup\limits_{\mathcal{N}_{p,q}(c,\mu) \rightarrow 0}s_{\mu,c} =1, $$
and hence $s_{\mu,c} \rightarrow 1$ as $\mathcal{N}_{p,q}(c,\mu) \rightarrow 0$. \qed

\begin{Lem}\label{Lemma5.5}
Let $\max\{ \mathcal{N}_{p,q}(c,\mu), \mathcal{M}_q(c,\mu)\} \rightarrow 0$. Then $v_{\mu,c} \rightarrow Z$ strongly in $H^1(\RN)$.
\end{Lem}
\Proof Adopt the similar argument of (\ref{LIMITG}), one has
\begin{equation}\label{LIMITZ}
    \int_{\RN}G_{\mu,c}(s_{\mu,c}Z)dx=\frac{1}{q}\|s_{\mu,c}Z\|_q^q+o_{\mu,c}(1).
\end{equation}
By (\ref{TheLowerConstant}), we obtain
\begin{equation}\label{LIMITVS}
     \Bigg|\int_{\RN}(I_\alpha * |v_{\mu,c}|^{\frac{N+\alpha}{N}})|v_{\mu,c}|^{\frac{N+\alpha}{N}}dx \Bigg| \leq S_{\alpha}^{-\frac{N+\alpha}{N}},
\end{equation}\begin{equation}\label{LIMITZS}
     \Bigg| \int_{\RN}(I_\alpha * |Z|^{\frac{N+\alpha}{N}})|Z|^{\frac{N+\alpha}{N}}dx \Bigg| \leq S_{\alpha}^{-\frac{N+\alpha}{N}}.
\end{equation}
It follows from (\ref{LIMITfunctional}), (\ref{Zfunctional}), (\ref{LIMITG}), (\ref{LIMITZ}), (\ref{LIMITVS}) and (\ref{LIMITZS}) that
$$ I(s_{\mu,c}Z)=J_{\mu,c}(s_{\mu,c}Z)+o_{\mu,c}(1) \quad \text{and}\quad J_{\mu,c}(t_{\mu,c}v_{\mu,c})=I(t_{\mu,c}v_{\mu,c})+o_{\mu,c}(1).  $$
By the similar argument of Lemma \ref{TheMaxFiber}, we get
$ I(Z) \geq I(s_{\mu,c}Z)$ and $ J_{\mu,c}(v_{\mu,v}) \geq J_{\mu,c}(t_{\mu,c}v_{\mu,c})$.
Thus, it follows that
\begin{equation*}
 \aligned e(1,1)+o_{\mu,c}(1) \geq I(t_{\mu,c}v_{\mu,c}) \geq  e(1,1).   \endaligned
\end{equation*}
Then $\{t_{\mu,c}v_{\mu,c}\} \subset \mathcal{M}$ is a bounded minimizing sequence of $e(1,1)$, since $t_{\mu,c}v_{\mu,c}(x)$ is radial symmetric, nonincreasing function of $|x|$ for any $\mu>0$
and $c>0$, it is standard to show $t_{\mu,c}v_{\mu,c}\rightarrow Z$ strongly in $H^1(\RN)$ by the uniqueness of the minimizer.
Since $t_{\mu,c}\rightarrow 1$, we get $v_{\mu,c} \rightarrow Z$ as $\mathcal{N}_{p,q}(c,\mu) \rightarrow 0$ and $ \mathcal{M}_q(c,\mu) \rightarrow 0$.
\qed

\noindent\textbf{Proof of Theorem \ref{ap1}.}\ From (\ref{Thefo}) we have
$$Z = \|U\|_2^{\frac{4}{N(q-2)-4}}U(\|U\|_2^{\frac{2(q-2)}{N(q-2)-4}}x)\quad \text{and}\quad \Lambda=\|U\|_2^{\frac{4(q-2)}{N(q-2)-4}}.  $$
Adopting the similar argument in Lemma \ref{aspG}, we have
$$ \int_{\RN}g_{\mu,c}(v_{\mu,c})v_{\mu,c}dx = \|v_{\mu,c}\|_q^q + o_{\mu,c}(1).   $$
and
$$\mu^{\frac{4}{N(q-2)-4}} c^{\frac{4(q-2)}{N(q-2)-4}}\lambda_{\mu,c} \rightarrow \|U\|_{2}^{\frac{4(q-2)}{N(q-2)-4}} $$
as $\max\{ \mathcal{N}_{p,q}(c,\mu), \mathcal{M}_q(c,\mu)\} \rightarrow 0$.
The proof is complete. \qed

Now let $(u_{\mu,c}, \lambda_{\mu,c})$ be a family of solutions obtained in Theorem \ref{theresult2} with $\gamma=1$, we introduce a new rescaling
\begin{equation}\label{Vscaling}
  w_{\mu,c}:= \mu^{\frac{N}{N(p-2)-4}} \lambda_{\mu,c}^{\frac{N}{\alpha}\frac{2}{N(p-2)-4}}u_{\mu,c}\left( \mu^{\frac{2}{N(p-2)-4}} \lambda_{\mu,c}^{\frac{N}{\alpha}\frac{p-2}{N(p-2)-4}} x \right).
\end{equation}
Then $w_{\mu,c}$ satisfies the equation
\begin{equation}\label{wequation}
  \aligned
 - \mu^{-\frac{4}{N(p-2)-4}} \lambda_{\mu,c}^{-\frac{N}{\alpha}\frac{2(p-2)}{N(p-2)-4}-1} \Delta w +  w
                                  = (I_{\alpha}*|w_{\mu,c}|^{\2})|w|^{\2-2}w + \tilde{g}_{\mu,c}(w) \endaligned
\end{equation}
where $\tilde{g}_{\mu,c}(w):=\mu^{\frac{N(p-1)-4}{N(p-2)-4}} \lambda_{\mu,c}^{\frac{N}{\alpha}\frac{2}{N(p-2)-4}-1 } g(\mu^{-\frac{N}{N(p-2)-4}} \lambda_{\mu,c}^{-\frac{N}{\alpha}\frac{2}{N(p-2)-4}} w)$. The corresponding action functional is defined by
$$\aligned \tilde{A}_{\mu,c}(w):= &\frac{1}{2}\mu^{-\frac{4}{N(p-2)-4}} \lambda_{\mu,c}^{-\frac{N}{\alpha}\frac{2(p-2)}{N(p-2)-4}} \|\nabla w\|_2^2 + \frac{1}{2}\|w\|_2^2 \\
                          & - \int_{\RN}\tilde{G}_{\mu,c}(w)dx - \frac{N}{2(N+\alpha)} \int_{\RN}(I_{\alpha}*|w|^{\2})|w|^{\2}dx.   \endaligned $$
where $\tilde{G}(w)= \mu^{\frac{Np-4}{N(p-2)-4}}  \lambda_{\mu,c}^{\frac{N}{\alpha}\frac{4}{N(p-2)-4}-1}G(\mu^{-\frac{N}{N(p-2)-4}} \lambda_{\mu,c}^{-\frac{N}{\alpha}\frac{2}{N(p-2)-4}} w). $
\begin{Lem}\label{uclargeequality}
For $\mathcal{M}_p(c,\mu)$ sufficiently large, there hold
\begin{enumerate}
  \item [(1)]$\|\nabla w_{\mu,c}\|_2^2 = \mu^{\frac{4}{N(p-2)-4}} \lambda_{\mu,c}^{\frac{N}{\alpha}\frac{2N-p(N-2)}{N(p-2)-4}}\|\nabla u_{\mu,c}\|_2^2 \sim 1$;
  \item [(2)]$\|w_{\mu,c} \|_2^2=\lambda_{\mu,c}^{-\frac{N}{\alpha}} \|u_{\mu,c} \|_2^2 \sim 1$;
  \item [(3)]$\int_{\RN}(I_{\alpha}*|w_{\mu,c}|^{\2})|w_{\mu,c}|^{\2}dx =\lambda_{\mu,c}^{-\frac{N+\alpha}{\alpha}} \int_{\RN}(I_{\alpha}*|u_{\mu,c}|^{\2})|u_{\mu,c}|^{\2}dx$;
  \item [(4)]$\int_{\RN}\tilde{G}(w_{\mu,c})dx = \mu \lambda_{\mu,c}^{-\frac{N+\alpha}{\alpha}}\int_{\RN}G(u_{\mu,c})dx$.
\end{enumerate}
\end{Lem}
\Proof It follows from direct calculation and (\ref{TheCharLag2}).  \qed

\begin{Lem}\label{Lemma5.7}
Let $\min\{ \mathcal{N}_{p,q}(c,\mu), \mathcal{M}_p(c,\mu) \}\rightarrow +\infty$, then
\begin{equation}\label{LIMITG2}
 \mu^{\frac{4}{N(p-2)-4}} \lambda_{\mu,c}^{\frac{N}{\alpha}\frac{2(p-2)}{N(p-2)-4}+1}\cdot \int_{\RN}\tilde{G}_{\mu,c}(w_{\mu,c})dx=\frac{1}{p}\|w_{\mu,c}\|_p^p+o_{\mu,c}(1)
\end{equation}
and
\begin{equation}\label{LIMITBarG2}
\aligned  &\mu^{\frac{4}{N(p-2)-4}} \lambda_{\mu,c}^{\frac{N}{\alpha}\frac{2(p-2)}{N(p-2)-4}+1} \cdot \frac{N}{2}\int_{\RN}\left[\tilde{g}_{\mu,c}(w_{\mu,c})w_{\mu,c}-2\tilde{G}_{\mu,c}(w_{\mu,c})\right]dx\\ =&\frac{N(p-2)}{2p}\|w_{\mu,c}\|_p^p+o_{\mu,c}(1) \endaligned
\end{equation}
where $o_{\mu,c}(1)\rightarrow 0$ as $\min\{ \mathcal{N}_{p,q}(c,\mu), \mathcal{M}_p(c,\mu) \}\rightarrow +\infty$.
\end{Lem}
\Proof Let $g(t):=t^{p-1}+\tilde{h}(t)$ for $t \geq 0$, then $G(t)=\frac{1}{p}t^p + \tilde{H}(t)$, and
$$ \lim_{t \rightarrow 0}\tilde{h}(t)/|t|^{p-1}=  \lim_{t \rightarrow 0}\tilde{H}(t)/|t|^p= 0. $$
Therefore, for any $\varepsilon>0$, there exists a constant $C_{\varepsilon}>0$ such that
$$ \Big|G(t)-\frac{1}{p}|t|^p \Big| \leq |\tilde{H}(t)| \leq \varepsilon|t|^{p} + C_{\varepsilon}|t|^{q}, $$
and hence we deduce that
\begin{equation*}
\aligned & \Bigg|\mu^{\frac{4}{N(p-2)-4}} \lambda_{\mu,c}^{\frac{N}{\alpha}\frac{2(p-2)}{N(p-2)-4}+1}\cdot \int_{\RN}\tilde{G}_{\mu,c}(w_{\mu,c})dx-\frac{1}{p}\|w_{\mu,c}\|_p^p \Bigg|  \\
\leq & \varepsilon\|w_{\mu,c}\|_p^p +\left[ \mu \lambda^{\frac{2}{\alpha}}\right]^{\frac{N(p-q)}{N(p-2)-4}}C_{\varepsilon}\|w_{\mu,c}\|_{q}^{q}. \endaligned
\end{equation*}
Similarly, for any $\varepsilon>0$, there exists a constant $C_{\varepsilon}>0$ such that
$$ \Big|\frac{N}{2}\bar{G}(t)-\frac{N(p-2)}{2p}|t|^p \Big| \leq \Big|\frac{N}{2}\tilde{h}(t)t-N\tilde{H}(t)\Big| \leq \varepsilon|t|^{p} + C_{\varepsilon}|t|^{q}, $$
and hence
$$\aligned & \mu^{\frac{4}{N(p-2)-4}} \lambda_{\mu,c}^{\frac{N}{\alpha}\frac{2(p-2)}{N(p-2)-4}+1} \cdot \frac{N}{2}\int_{\RN}\left[\tilde{g}_{\mu,c}(w_{\mu,c})w_{\mu,c}-2\tilde{G}_{\mu,c}(w_{\mu,c})\right]dx-\frac{N(p-2)}{2p} \|w_{\mu,c}\|_p^p \\
     \leq &\varepsilon \|w_{\mu,c}\|_p^p +\left[ \mu \lambda^{\frac{2}{\alpha}}\right]^{\frac{N(p-q)}{N(p-2)-4}}C_{\varepsilon}\|w_{\mu,c}\|_{q}^{q}. \endaligned$$
Since $\left[ \mu \lambda^{\frac{2}{\alpha}}\right]^{\frac{N(p-q)}{N(p-2)-4}} \rightarrow 0$ as $\min\{ \mathcal{N}_{p,q}(c,\mu), \mathcal{M}_p(c,\mu) \}\rightarrow +\infty$, the rest of the proof is similar to Lemma \ref{aspG} and we omit it. \qed

\begin{Lem}\label{Lemma5.8}
Let $\min\{ \mathcal{N}_{p,q}(c,\mu), \mathcal{M}_p(c,\mu) \}\rightarrow +\infty$, then $w_{\mu,c} \rightarrow V_{\delta_0} $ strongly in $H^1(\RN)$ where $V_{\delta_0}$ is a minimizer of $S_{\alpha}$ such that $\|V_{\delta_0}\|_2^2 =S_{\alpha}^{\frac{N+\alpha}{\alpha}}$ and $\|\nabla V_{\delta_0} \|_2^2=\frac{N(p-2)}{2p}\| V_{\delta_0}\|_p^p$.
\end{Lem}
\Proof From Lemma \ref{uclargeequality} and (\ref{LIMITG2}), we deduce
$$ \aligned \tilde{A}_{\mu,c}(w_{\mu,c}) = &\frac{1}{2}\mu^{-\frac{4}{N(p-2)-4}} \lambda_{\mu,c}^{-\frac{N}{\alpha}\frac{2(p-2)}{N(p-2)-4}} \|\nabla w_{\mu,c}\|_2^2 + \frac{1}{2}\|w_{\mu,c}\|_2^2 \\
                          & + \int_{\RN}\tilde{G}_{\mu,c}(w_{\mu,c})dx + \frac{N}{2(N+\alpha)} \int_{\RN}(I_{\alpha}*|w_{\mu,c}|^{\2})|w_{\mu,c}|^{\2}dx \\
                                         =& \frac{1}{2}\|w_{\mu,c}\|_2^2 -\frac{N}{2(N+\alpha)} \int_{\RN}(I_{\alpha}*|w_{\mu,c}|^{\2})|w_{\mu,c}|^{\2}dx + o_{\mu,c}(1),  \endaligned $$
and
$$  \tilde{A}_{\mu,c}(w_{\mu,c}) = \lambda_{\mu,c}^{-\frac{N+\alpha}{\alpha}} \left(E_{\mu,c}(u_{\mu,c}) +\frac{1}{2}\lambda_{\mu,c}\| u_{\mu,c}\|_2^2 \right)=\frac{\alpha}{2(N+\alpha)}S_{\alpha}^{\frac{N+\alpha}{\alpha}}+o_{\mu,c}(1).  $$
On the other hand, since $w_{\mu,c}$ solves (\ref{wequation}) and for $\min\{ \mathcal{N}_{p,q}(c,\mu), \mathcal{M}_p(c,\mu) \}$ sufficiently large, $\{w_{\mu,c}\}$ is bounded in $H^1(\RN)$, we deduce that
$$  \aligned   0= \tilde{A}_{\mu,c}'(w_{\mu,c})= \frac{1}{2}w_{\mu,c} - (I_{\alpha}*|w_{\mu,c}|^{\2})|w_{\mu,c}|^{\2-2}w_{\mu,c}+o_{\mu,c}(1) \quad \text{in\ }\left(H^1(\RN)\right)^{*}.   \endaligned  $$
Let us define a functional $\tilde{A}_0$ on $H^1(\RN)$ by $$\tilde{A}_0(w)=\frac{1}{2}\| w\|_2^2 - \frac{N}{2(N+\alpha)}\int_{\RN}(I_{\alpha}*|w|^{\2})|w|^{\2}dx, $$
then $\{w_{\mu,c}\}$ is a bounded Palais-Smale sequence for $\tilde{A}_0$ at the level $M_0=\frac{\alpha}{2(N+\alpha)}S_{\alpha}^{\frac{N+\alpha}{\alpha}}$. We may assume $w_{\mu,c} \rightarrow w_0$, then it follows from Lemma \ref{TheCompactnessTool} that $w_{\mu,c} \rightarrow w_0$ in $L^s(\RN)$ for any $ s \in (2,2^*)$.
By Lemma \ref{uclargeequality} and (\ref{LIMITBarG2}), we deduce
$$  \aligned &\|\nabla w_{\mu,c}\|_2^2 = \mu^{\frac{4}{N(p-2)-4}} \lambda_{\mu,c}^{\frac{N}{\alpha}\frac{2(p-2)}{N(p-2)-4}+1} \cdot \frac{N}{2}\int_{\RN}\left[\tilde{g}_{\mu,c}(w_{\mu,c})w_{\mu,c}-2\tilde{G}_{\mu,c}(w_{\mu,c})\right]dx  \\
    =&\frac{N(p-2)}{2p}\|w_{\mu,c}\|_p^p+o_{\mu,c}(1)   \endaligned $$
and hence $w_0 \not\equiv 0$.
Then by the weakly convergent and a standard argument it follows that $\tilde{A}_0'(w_0)=0$ and $\tilde{A}_0(w_0) \geq M_0$.
Arguing as in the proof of \cite[Lemma 4.8 and Theorem 2.1]{MaMoroz}, we obtain that there exists a family $\{ \zeta_{\mu,c} \} \subset \mathbb{R}$ such that
  $$ \zeta_{\mu,c} ^{\frac{N}{2}}u_{\mu,c}( \zeta_{\mu,c}  x) \rightarrow V_{\delta_0}  $$
as $\min\{ \mathcal{N}_{p,q}(c,\mu), \mathcal{M}_p(c,\mu)\}\rightarrow +\infty$ where
$$ \zeta_{\mu,c} \sim \mu^{\frac{2}{N(p-2)-4}}\lambda_{\mu,c}^{\frac{1}{2}+\frac{N}{\alpha}\frac{p-2}{N(p-2)-4}} \sim \mathcal{M}_p(c,\mu)^{\frac{1}{2}} $$
and $V_{\delta_0}$ is a minimizer of $S_\alpha$ such that $\|V_{\delta_0}\|_2^2 =S_{\alpha}^{\frac{N+\alpha}{\alpha}}$ and $\|\nabla V_{\delta_0} \|_2^2=\frac{N(p-2)}{2p}\| V_{\delta_0}\|_p^p$.
\qed

\noindent\textbf{Proof of Theorem \ref{ap2}.} It follows from Lemma \ref{Lemma5.8}.  \qed

\section{On the power type nonlinearity}

Now we consider $g(u)=\mu |u|^{q-2}u$ with $ 2+\frac{4}{N}<q< 2^*$. In this case, we deduce
$$ E_{\mu,c}(u)=\frac{1}{2}\|\nabla u\|_2^2 - \frac{\mu}{q}\|u \|_q^q  - \frac{N\gamma}{2(N+\alpha)}\II, $$
and
$$ \mathcal{P}(\mu,c)=\left\{ u\in S(c): \|\nabla u\|_2^2 = \mu \frac{N(q-2)}{2q}\|u \|_q^q  \right\}. $$

\noindent\textbf{Proof of Theorem \ref{Homoresult}.} (1) Since $g$ satisfies $(G1)$--$(G5)$, then the proof is similar to the previous argument.

(2) Since $\tilde{f} \equiv 0$, then we have $t_{\mu,c}=1$ in Lemma \ref{Lemma5.3} and $s_{\mu,c}=1$ in Lemma \ref{Lemma5.4}. Hence, we only need $\mathcal{M}_q(c,\mu)\rightarrow 0$ in Lemma \ref{Lemma5.5}.
On the other hand, since $\tilde{h} \equiv 0$ we only need $\mathcal{M}_q(c,\mu)\rightarrow +\infty$ Lemma \ref{Lemma5.7} and Lemma \ref{Lemma5.8}. The rest of the proof is similar to the previous argument. \qed
\begin{Rem}
In fact, through some tedious calculations, one can ascertain the explicit formulations for $\omega_1, \omega_2$. However, it is noteworthy that
$\omega_1, \omega_2$ may not necessarily serve as the threshold parameters of $\mathcal{M}_q(c,\mu)$ that determine the existence of the normalized ground state solution.
In practice, the assumption $\mathcal{M}_q(c,\mu) \not \in [\omega_1,\omega_2]$ is a technical condition to ensure
$$  \mathcal{E}(\mu,c)+\frac{1}{2}\lambda_{\mu,c}c^2 < a(\mu,c)+\frac{\alpha}{2(N+\alpha)}\gamma^{-\frac{N}{\alpha}}(\lambda_{\mu,c} S_{\alpha})^{\frac{N+\alpha}{\alpha}}. $$

\noindent Consequently, the scenario where $\mathcal{M}_q(c,\mu) \in [\omega_1,\omega_2]$ remains entirely unexplored within the current analytical framework.
We believe it would be interesting to investigate in this direction, and we propose the following conjecture:

\noindent\textbf{Conjecture.} Let $N \geq 1$, $\gamma > 0$, $\mu>0$, $c>0$, and $2+ \frac{4}{N} < q <2^* $. Then there exist two positive constants $ \kappa_1, \kappa_2$ independent of $\mu$, $c$ such that $\omega_1 \leq \kappa_1 \leq \kappa_2 \leq \omega_2$ and
\begin{enumerate}
  \item [(1)] Equ. $(\mathcal{A}_{\mu,c})$ admits a normalized ground state solution if $\mathcal{M}_q(c,\mu) \not \in [\kappa_1,\kappa_2]$;
  \item [(2)] Equ. $(\mathcal{A}_{\mu,c})$ admits no normalized ground state solution if $\mathcal{M}_q(c,\mu) \in [\kappa_1,\kappa_2]$.
\end{enumerate}
\end{Rem}
Then we study the non-existence of least action solutions to $(B_{\eta})$
\begin{equation*}
    \left\{ \aligned  &-\Delta u + u = \eta |u|^{q-2}u + (I_\alpha * |u|^{\frac{N+\alpha}{N}})|u|^{\frac{N+\alpha}{N}-2}u   & \text{in\ \ }  \mathbb{R}^N,\\
    & u \in H^1(\mathbb{R}^N) \endaligned \right.
    \end{equation*}
for $\eta >0$. The corresponding action functional is
\begin{equation}\label{Theactionfunctional1}
  A_{\eta}(u):=\frac{1}{2}\|\nabla u\|_2^2 + \frac{1}{2}\| u \|_2^2 - \frac{\eta}{q}\| u\|_q^q - \frac{N}{2(N+\alpha)}\II,
\end{equation}
and by \cite{Tangxianhua1}, the least action level $a(\eta)$ can be characterized as
$$ a(\eta)= \inf_{u \in \mathcal{N}(\eta)} A_{\eta}(u) $$
where
$$ \mathcal{N}(\eta):=\left\{ u \in H^1(\RN)\setminus\{0\}:\|\nabla u\|_2^2 + \|u \|_2^2 = \eta \|u\|_q^q + \II  \right\} $$
is the Nehari manifold.

By using the associated fibering map, it is easy to show that
\begin{Lem}\label{fiberfixedfrequency}
Let $N \geq 1$, $\eta \geq 0$, $q \in [2+\frac{4}{N},2^*)$. Then
\begin{enumerate}
  \item [(1)] for any $u \in H^1(\RN)\setminus\{0\}$, there exists unique $t_{\eta,u}$ such that $t_{\eta,u}u\in\mathcal{N}(\eta)$, and
  $$ a(\eta)=\inf_{u\in H^1(\RN)\setminus\{0\}}\max_{t > 0}A_{\eta}(tu)>0;   $$
  \item [(2)] $a(0)=\frac{\alpha}{2(N+\alpha)}S_{\alpha}^{\frac{N+\alpha}{\alpha}}$, $a(\eta) \rightarrow 0^+$ as $\eta \rightarrow +\infty$;
  \item [(3)] the mapping $\eta \mapsto a(\eta)$ is nonincreasing and continuous on $[0,+\infty)$;
  \item[(4)] $(\mathcal{B}_{\eta})$ admits a positive, radially symmetric, radially nonincreasing least action solution if $a(\eta)< \frac{\alpha}{2(N+\alpha)}S_{\alpha}^{\frac{N+\alpha}{\alpha}}$.
\end{enumerate}
\end{Lem}
\begin{Lem}\label{eta1def}
Let $N \geq 1$, $q \in [2+\frac{4}{N}, 2^*)$. Then $a(\eta)=\frac{\alpha}{2(N+\alpha)}S_{\alpha}^{\frac{N+\alpha}{\alpha}}$ for $\eta > 0$ sufficiently small.
\end{Lem}
\Proof We argue in the contrary by supposing that there exists $\eta_n \rightarrow 0$ as $n \rightarrow \infty$ such that $a(\eta_n)< \frac{\alpha}{2(N+\alpha)}S_{\alpha}^{\frac{N+\alpha}{\alpha}}$.
Then it follows from Lemma \ref{fiberfixedfrequency} (4) that $a(\eta)$ is attainted by a least action solution of $(B_{\eta_n})$. We denote this solution by $w_n$, and we deduce
\begin{equation}\label{thecont1}
\|\nabla w_n\|_2^2 + \|w_n \|_2^2 =  \eta_n \| w_n\|_q^q + \int_{\RN}(I_{\alpha}*|w_n|^{\2})|w_n|^\2 dx,
\end{equation}
\begin{equation}\label{thecont2}
  \|\nabla w_n\|_2^2 = \eta_n \frac{N(q-2)}{2q}\|w_n \|_q^q.
\end{equation}
It follows from (\ref{Theactionfunctional1}) and (\ref{thecont1}) that
$$ a(\eta_n) \geq \left(\frac{1}{2}- \max\left\{\frac{1}{q},\frac{N}{2(N+\alpha)}\right\} \right)\left(\|\nabla w_n\|_2^2 + \|w_n\|_2^2 \right), $$
and hence, $\{w_n \}$ is bounded in $H^1(\RN)$. As a consequence, by  the Gagliardo-Nirenberg inequality  (Lemma \ref{GNS}) and (\ref{thecont2}) one has
 \begin{equation*}
\aligned 1 & \leq \eta_n \frac{N(q-2)}{2q} C_{N,q}^q\| w_n \|_2^{\frac{2q-N(q-2)}{2}} \| \nabla w_n\|_2^{\frac{N(q-2)-2}{2}} \endaligned
\end{equation*}
which is a contradiction with $\eta_n \rightarrow 0$ as $n \rightarrow \infty$.
Then it follows from Lemma \ref{fiberfixedfrequency} (2)(3) that $a(\eta) = \frac{\alpha}{2(N+\alpha)}S_{\alpha}^{\frac{N+\alpha}{\alpha}}$
for $\eta>0$ sufficiently small. \qed

As a consequence of Lemma \ref{fiberfixedfrequency} (2) and Lemma \ref{eta1def}, we obtain
$$ 0<\eta_1: =\sup\left\{\eta>0 : a(\eta) = \frac{\alpha}{2(N+\alpha)}S_{\alpha}^{\frac{N+\alpha}{\alpha}}   \right\} < +\infty. $$

\begin{Lem}\label{Nonexistence}
Let $N \geq 1$, $q \in [2+\frac{4}{N}, 2^*)$. Then
\begin{enumerate}
  \item [(1)] $(\mathcal{B}_{\eta})$ admits no least action solution for $\eta < \eta_1$ in the case of $q \in [2+\frac{4}{N}, 2^*)$;
  \item [(2)] $(\mathcal{B}_{\eta})$ admits a least action solution for $\eta \geq \eta_1$ in the case of $q \in (2+\frac{4}{N}, 2^*)$.
\end{enumerate}
\end{Lem}
\Proof (1) Suppose the contrary that
$(\mathcal{B}_{\eta})$ admits a least action solution $w$ for some $0 < \eta <\eta_1$. Then for any $\eta^*>\eta$, there exists $t^*>0$ such that $t^*w \in \mathcal{N}(\eta^*)$, and this implies
\begin{equation}\label{supdef}
\aligned  a(\eta)=A_{\eta}(w)\geq & A_{\eta}(t^*w) = A_{\eta^*}(t^*w)+(\eta^*-\eta)\frac{(t^*)^q}{q}\|w \|_q^q \\
                       > &A_{\eta^*}(t^*w) \geq a(\eta^*)\endaligned
\end{equation}
from Lemma \ref{fiberfixedfrequency} (1). Since $ a(\eta)=\frac{\alpha}{2(N+\alpha)}S_{\alpha}^{\frac{N+\alpha}{\alpha}}$ for $\eta < \eta_1$, we obtain (\ref{supdef}) is a contradiction with the definition of $\eta_1$.
Thus, there is no least action solution of $(\mathcal{B}_{\eta})$ for $\eta < \eta_1$.

(2) Let $\{ \tilde{\eta}_n \}$ be a sequence such that $\tilde{\eta}_n \rightarrow \eta_1^+$ as $n \rightarrow \infty$ and
$\{w_n\}$ be a sequence of least
action solutions of $(\mathcal{B}_{\tilde{\eta}_n})$, then it is standard to show $\{w_n\}$ is a bounded Palais-Smale sequence of $A_{\eta_1}$ at the
level $a(\eta_1)$ from Lemma \ref{fiberfixedfrequency} (3). Moreover, since $w_n$ is a solution to $(\mathcal{B}_{\tilde{\eta}_n})$, one has
\begin{equation}\label{Theupper}
  \|\nabla w_n \|_2^2 = \tilde{\eta_n} \frac{N(q-2)}{2q}\|w_n \|_q^q,
\end{equation}
and
\begin{equation}\label{Thelower}
 \|\nabla w_n\|_2^2 +\|w_n \|_2^2 \leq (\eta_1+1) \|w_n \|_q^q + \int_{\RN}(I_{\alpha}*|w_n|^{\2})|w_n|^{\2}dx.
\end{equation}
Now we may assume $w_n \rightharpoonup w_0 $ in $H^1(\RN)$, by Lemma \ref{TheCompactnessTool}, one has $w_n \rightarrow w_0$ in $L^q(\RN)$ and $w_n \rightharpoonup w_0$ a.e in $\RN$. From (\ref{Thelower}), the Hardy-Littlewood-Sobolev
inequality and Sobolev embedding theorem, we deduce
\begin{equation}\label{thelowerL2}
  \|\nabla w_n \|_2^2 + \|w_n \|_2^2 \gtrsim 1 \quad \text{as\ } n \rightarrow \infty.
\end{equation}
If $w_0=0$, since $w_n$ is radial and nonincreasing function of $|x|$ for every $n \geq 0$, then it follows from Lemma \ref{TheCompactnessTool}, (\ref{Theupper}) and (\ref{thelowerL2}) that
\begin{equation}\label{imposs}
  \|\nabla w_n \|_2^2 = \|w_n \|_q^q =o(1), \quad \|w_n \|_2^2 \gtrsim 1 \quad \text{as\ } n \rightarrow \infty.
\end{equation}
However, by the Gagliardo-Nirenberg inequality (Lemma \ref{GNS}), one has
$$ \|\nabla w_n\|_2^2 \lesssim \|w_n \|_q^q \lesssim \|\nabla w_n\|_2^{\frac{N(q-2)}{2}} \quad \text{as\ } n \rightarrow \infty $$
which implies $\|\nabla w_n \|_2^2 \gtrsim 1$ as $n \rightarrow \infty$ and contradicts (\ref{imposs}). Thus, one has $w_0 \neq 0$ and $w_0$ is a nontrivial positive solution which is
radially symmetric and radially nonincreasing. Besides, we obtain
$$\aligned a(\eta_1)+o(1) & =\left(\frac{1}{2}- \max\left\{\frac{1}{q},\frac{N}{2(N+\alpha)}\right\} \right)\left(\|\nabla w_n\|_2^2 + \|w_n\|_2^2 \right)\\
                              &\geq \left(\frac{1}{2}- \max\left\{\frac{1}{q},\frac{N}{2(N+\alpha)}\right\} \right)\left(\|\nabla w_0\|_2^2 + \|w_0\|_2^2 \right)\\
                              &=A_{\eta_1}(w_0)\geq a(\eta_1). \endaligned$$
Then $w_0$ is a least action solution of $(\mathcal{B}_{\eta_1})$ and this completes the proof. \qed

Now we explore the multiplicity of positive solutions to $(\mathcal{B}_{\eta})$.
Let $\{(u_{\mu}, \lambda_{\mu})\}$ denotes the family of the normalized ground state solutions to $\left(\mathcal{A}_{\mu,c}\right)$ obtained in Theorem \ref{Homoresult} for $c>0$ fixed, $\gamma=1$ and $\mu \rightarrow +\infty$,
we introduce another rescaling by
$$ \bar{u}_{\mu}:= \lambda_{\mu}^{-\frac{N(2+\alpha)}{4\alpha}}u_{\mu}(\lambda_{\mu}^{-\frac{1}{2}}x).      $$
Then $\bar{u}_{\mu}$ satisfies the equation $(\mathcal{B}_{\eta_{\mu}})$
$$ \aligned &-\Delta u + u = \eta_{\mu} |u|^{q-2}u + (I_{\alpha}*|u|^{\2})|u|^{\2-2}u &\text{in\ }\RN
\endaligned  $$
where $\eta_{\mu}:=\mu \lambda_{\mu}^{\bar{q}}$ and $\bar{q}=\frac{N(q-2)-4}{4}+\frac{N(q-2)}{2\alpha}$. Besides, it is easy to see that $\eta_{\mu} \rightarrow +\infty$ as $\mu \rightarrow +infty$.
By direct calculations, we have
\begin{Lem}\label{6.1}
\begin{enumerate}
  \item [(1)] $\|\nabla \bar{u}_{\mu}  \|_2^2 = \lambda_{\mu}^{-\frac{N+\alpha}{\alpha}} \|\nabla u_{\mu}  \|_2^2$, $\|\bar{u}_{\mu}\|_2^2=\lambda_{\mu}^{-\frac{N}{\alpha}}\|u_{\mu}  \|_2^2$;
  \item [(2)] $  \|\bar{u}_{\mu}\|_q^q=\lambda_{\mu}^{-\frac{N(2+\alpha)q}{4\alpha}+\frac{N}{2}}\|u_{\mu}  \|_q^q $;
  \item [(3)] $ \int_{\RN}(I_{\alpha}*|\bar{u}_{\mu}|^{\2})|\bar{u}_{\mu}|^{\2}dx = \lambda_{\mu}^{-\frac{N+\alpha}{\alpha}} \int_{\RN}(I_{\alpha}*|u_{\mu}|^{\2})|u_{\mu}|^{\2}dx $.
\end{enumerate}
\end{Lem}

\noindent\textbf{Proof of Theorem \ref{Mutli}.} (1) (2) follow from Lemma \ref{Nonexistence} and \ref{fiberfixedfrequency} (2)(4).

(3) From Lemma \ref{6.1}, (\ref{TheCharLag2}) and (\ref{TheCharOfEnergy2}), we deduce
\begin{equation}
\aligned A_{\eta_{\mu}}(\bar{u}_{\mu})=& \frac{1}{2}\lambda_{\mu}^{-\frac{N+\alpha}{\alpha}} \|\nabla u_{\mu}  \|_2^2 + \frac{1}{2}\lambda_{\mu}^{-\frac{N}{\alpha}}\|u_{\mu}  \|_2^2 \\
                          &-\frac{\mu}{q}\lambda_{\mu}^{\bar{q}}\lambda_{\mu}^{-\frac{N(2+\alpha)q}{4\alpha}+\frac{N}{2}}\|u_{\mu}  \|_q^q-\frac{N}{2(N+\alpha)}\lambda_{\mu}^{-\frac{N+\alpha}{\alpha}} \int_{\RN}(I_{\alpha}*|u_{\mu}|^{\2})|u_{\mu}|^{\2}dx\\
                         =& \frac{1}{2}\lambda_{\mu}^{-\frac{N+\alpha}{\alpha}} \|\nabla u_{\mu}  \|_2^2 + \frac{1}{2}\lambda_{\mu}^{-\frac{N}{\alpha}}\|u_{\mu}  \|_2^2 \\
                           &-\frac{\mu}{q}\lambda_{\mu}^{-\frac{N+\alpha}{\alpha}} \|u_{\mu}  \|_q^q-\frac{N}{2(N+\alpha)}\lambda_{\mu}^{-\frac{N+\alpha}{\alpha}} \int_{\RN}(I_{\alpha}*|u_{\mu}|^{\2})|u_{\mu}|^{\2}dx\\
                          =&\lambda_{\mu}^{-\frac{N+\alpha}{\alpha}}\left( E_{\mu,c}(u_\mu)+\frac{1}{2}\lambda_{\mu}c^2  \right) \\
                         \rightarrow & \frac{\alpha}{2(N+\alpha)}S_{\alpha}^{\frac{N+\alpha}{\alpha}}\quad \text{as\ }\mu \rightarrow +\infty.\endaligned
                           \end{equation}

Since $\eta_{\mu} \rightarrow+\infty$ as $\mu\rightarrow +\infty$, then $(\mathcal{B}_{\eta})$ admits a positive radial solution $u_{\eta}:=\bar{u}_{\mu}$ for $\eta=\eta_{\mu}$ sufficiently large. Moreover, $(\mathcal{B}_{\eta})$ admits a least action solution $v_{\eta}$ for $ \eta $ sufficiently large such that $A_{\eta}(v_{\eta}) \rightarrow 0$ as $\eta \rightarrow +\infty$. Hence,  there exists $\eta_2$ such that $A_{\eta}(u_{\eta}) > A_{\eta}(v_{\eta})$ if $\eta > \eta_2$. This completes the proof.
                           \qed

\bibliographystyle{amsplain}

\end{document}